\documentclass{elsarticle}

\usepackage{lineno}
\usepackage{hyperref}
\modulolinenumbers[5]

\journal{Computer Methods in Applied Mechanics and Engineering}

\usepackage{amssymb}
\usepackage{amsmath,amssymb,amsfonts}
\usepackage[margin=3.5cm]{geometry}
\usepackage{lipsum}

\usepackage{academicons}
\usepackage{url}            
\usepackage{booktabs}      
\usepackage{amsfonts}       
\usepackage{nicefrac}       
\usepackage{microtype}      
\usepackage{lipsum}
\usepackage{graphicx}
\usepackage{subfig}
\usepackage[export]{adjustbox}
\usepackage{caption}
\usepackage{physics}
\usepackage{accents}
\usepackage{xcolor}

\newcommand{\RN}[1]{%
  \textup{\lowercase\expandafter{\romannumeral#1}}%
}
\usepackage{tikz}
\usepackage{scalerel}
\DeclareRobustCommand{\blackline}{\raisebox{2pt}{\tikz{\draw[-,black,solid,line width = 1pt](0,0) -- (5mm,0)}}}
\DeclareRobustCommand{\redline}{\raisebox{2pt}{\tikz{\draw[-,red,dashed,line width = 1pt](0,0) -- (5mm,0)}}}

\DeclareRobustCommand{\blueline}{\raisebox{2pt}{\tikz{\draw[-,blue,solid,line width = 1pt](0,0) -- (5mm,0)}}}
\definecolor{post_b}{rgb}{0.30,0.75,0.93}

\newcommand{\orcid}[1]{\href{https://orcid.org/#1}{\includegraphics[width=10pt]{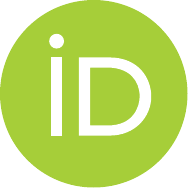}}}

\begin{document}

\begin{frontmatter}

\title{A Physics Informed Neural Network for Time--Dependent Nonlinear and Higher Order Partial Differential Equations}

\author[label2]{Revanth Mattey \orcid{0000-0002-6732-622X}}
\author[label2,label1]{Susanta Ghosh \orcid{0000-0002-6262-4121} \corref{cor1}}  
\cortext[cor1]{Corresponding author; Email: susantag@mtu.edu}

\address[label2]{Department of Mechanical Engineering--Engineering Mechanics, Michigan Technological University, MI, USA}
\address[label1]{The Center for Data Sciences, Michigan Technological University, MI, USA}

\begin{abstract}
A physics informed neural network (PINN) incorporates the physics of a system by satisfying its boundary value problem through a neural network's loss function. 
The PINN approach has shown great success in approximating the map between the solution of a partial differential equation (PDE) and its  spatio-temporal input. However, for strongly non-linear and higher order partial differential equations PINN's accuracy reduces significantly. To resolve this problem, we propose a novel PINN scheme that solves the PDE sequentially over successive time segments using a single neural network. The key idea is to re-train the same neural network for solving the PDE over successive time segments while satisfying the already obtained solution for all previous time segments. 
Thus it is named as backward compatible PINN (bc-PINN). To illustrate the advantages of bc-PINN, we have used the Cahn Hilliard and Allen Cahn equations, which are widely used to describe phase separation and reaction diffusion systems. Our results show significant improvement in accuracy over the PINN method while using a smaller number of collocation points. Additionally, we have shown that using the phase space technique for a higher order PDE could further improve the accuracy and efficiency of the bc-PINN scheme.
\end{abstract}


\end{frontmatter}

\section{Introduction}\label{sec:intro}
Traditional physics-based numerical methods for solving partial differential equations (PDEs) have found remarkable success in solving various science and engineering problems. These methods are accurate but computationally expensive for complex problems such as nonlinear PDEs and requires problem--specific techniques. In the last decade data driven methods have gained a lot of attention in almost all areas of science and engineering. Data driven methods for PDEs can help in identifying highly non-linear mappings (between the inputs and outputs) which can substitute or augment expensive physics based simulations. Due to their versatility and fast evaluation capabilities,  machine learning based models can be used as PDE solvers in situations when there is a requirement of large number of simulations such as the inverse--problem and homogenization \cite{MAO2020112789, Arbabi2020}. \\ 

Several data--driven techniques have been attempted to solve PDEs. For instance, the  Gaussian Process based approaches described in \cite{Raissi2017a,Raissi2017,Raissi2018,ATKINSON2019166,BILIONIS2013212}. Despite the ease of training for Gaussian Process, this approach did not gain as much popularity as neural network for solving PDEs due to its difficulties in handling high dimensional problems. \\ 

Among different data--driven techniques for PDEs the Physics Informed Neural Networks (PINN) has shown remarkable promise and versatility. PINN is a new class of  machine learning technique where a neural network's loss function is designed to satisfy the Initial Boundary Value Problem (IBVP) \cite{Raissi2019}. A PINN ``learns'' the non linear map between the  spatio--temporal input and the solution of the PDE in a given domain. PINN utilizes the automatic-differentiation capability \cite{baydin2018automatic} to compute the derivatives of the field variables.\\

Different variants of PINN  are shown to work effectively in solving many forward and inverse problems \cite{ZHANG2019108850,doi:10.1137/18M1225409,MENG2020109020}. Recently in \cite{Jagtap2020}, PINNs have been extended to satisfy various conservation laws while solving the PDEs. This approach is named as cPINNs. cPINNs solve the problem over several sub-domains and ensure flux continuity at the boundaries of the sub-domains.
While most of the PINN approaches solves the strong form of a PDE, it can also be used to solve the weak (variational) form of a PDE. Since the weak form incorporates the natural boundary conditions, the neural network solution only needs to satisfy the essential boundary conditions \textit{a priori}. This aspect is used in several numerical methods for PDEs such as finite element method. Due to this advantage of weak form over strong form the application of PINNs on the weak form have been in investigated in \cite{Karumuri2020}. In this study the authors have considered the variational form for stochastic PDEs and  applied the idea of PINNs to obtain the solution of the PDE. Also the corresponding uncertainty propagation through their model is presented in  \cite{Karumuri2020,Tripathy2018}. Uncertainty quantification provides the variation associated with the prediction of the model. It is particularly useful for systems where there is a high cost of data acquisition or lack of high resolution data \cite{Yang2019a}. 
The authors in \cite{yang2020bpinns} proposed a Bayesian approach for physics informed neural network to solve forward and inverse problems. \\ 

The promise and versatility of PINN have been demonstrated through its application for a wide range of problems. PINN has been used in  modelling subsurface transport phenomena \cite{HE2020103610}, approximating Euler equations for high speed flows \cite{MAO2020112789}, constitutive modeling of stress strain behavior in biological tissues \cite{LIU2020113402}, predicting arterial blood pressure from noisy MRI data of flow velocity \cite{KISSAS2020112623}, and cardiac activation mapping for diagnosing atrial fibrillation \cite{10.3389/fphy.2020.00042}. \\

In the present work, we demonstrate that the accuracy of the state-of-the-art PINN \cite{Raissi2019} suffers in the presence of  (\RN{1}) Strong Non-linearity, and (\RN{2}) Higher order partial differential operators. In order to illustrate the above, we chose the Allen Cahn equation having \textit{strong non-linearity} and Cahn Hilliard equation having \textit{strong non-linearity and fourth order derivative}. 
These are the two most widely used PDEs to study diffusion separation and multi phase flows \cite{doi:10.1142/S0218202511500138, deckelnick_dziuk_elliott_2005, doi:10.1098/rspa.1998.0273, CiCP-6-433}. To overcome the drawbacks of state of the art PINN, we have proposed an extension, which is named as backward compatible PINN (bc-PINN). The proposed bc-PINN solves the PDE over successive time segments by re-training the same neural network, where the key idea is: 
\begin{itemize}
    \item[] \hspace{-20pt} \emph{To ensure that the neural network can reproduce the solution for all the prior time segments while solving the PDE for a particular time segment.}
\end{itemize}
Henceforth, this idea is referred as \emph{backward compatibility}. Some of the main advantages of the proposed bc-PINN method are as follows: 
\begin{enumerate}
    \item It works for higher order and strongly nonlinear PDEs by using less number of iterations and collocation points while achieving significantly higher accuracy when compared to standard PINN.
    \item A single neural network is used for the entire domain and  continuity across the time segments is ensured for the predicted solution and its derivatives. 
\end{enumerate}

The rest of the paper is organized as follows: in section~(\ref{sec:DNN_PDE}) the PINN method is briefly reviewed; in section~(\ref{sec:BC-PINN}) the proposed bc-PINN method is presented;  in section~(\ref{sec:AC}) and ~(\ref{sec:CH}) the bc-PINN method is analyzed and compared against the PINN method for the Allen Cahn equation and Cahn Hilliard equation respectively. Finally, the conclusions are presented in section~(\ref{sec:conclusion}).


\section{A brief review of physics informed neural network (PINN) for partial differential equations}\label{sec:DNN_PDE}
Physics informed neural network (PINN) is a class of machine learning model where the governing PDE  is satisfied through the loss function of the neural network. The efficient optimization and prediction capabilities of neural network are exploited in the PINN approach. In PINN a neural network is trained to predict the solution at any point in the entire spatial--temporal domain.
Let's consider the general form of a $m^{\mathrm{th}}$ order partial differential equation (PDE):
\begin{equation}
    \centering
    h_t = F(h(x,t), h_x^{(1)}(x,t)\,, h_{x}^{(2)}(x,t),\cdots,h_x^{(m)}(x,t))\,, \;
    x\in\Omega \subset \mathbb{R}\,, \; 
    t \, \in \, (0,T]
    \label{eq:1}
\end{equation}
\sloppy Here, $\Omega$ is an open set of $\mathbb{R}$. \textit{F} is a non linear function of the solution $h(x,t)$ and it's  spatial derivatives ($h_x^{(1)}(x,t), h_x^{(2)}(x,t),\cdots,h_x^{(m)}(x,t)$) where $x$ and $t$ are the space and time coordinates respectively. 
The corresponding boundary conditions and initial conditions are
\begin{align}
    \begin{gathered}
    h(x,0) = \phi(x), \quad x \in \Omega \\
    h(-x,t) = h(x,t), \quad (x,t) \in \Gamma \times (0,T]  \\
    h_x^{(1)}(-x,t) = h_x^{(1)}(x,t)\, , \quad (x,t) \in \Gamma \times (0,T] 
    \end{gathered}
    \label{eq:2}
\end{align}
Where, $\Gamma$ is the boundary of $\Omega$. 
The PDE, the Initial and the Boundary Conditions (given by equation (\ref{eq:1}--\ref{eq:2})) form a  initial--boundary value problem (IBVP)  considered in this study. The boundary conditions are taken as periodic and the initial condition is a real function. \\

PINN  approximates the map between points in the  spatio-temporal domain to the solution of the PDE. The parameters of the  neural network are randomly initialized and iteratively updated by minimizing the loss function that enforces the PDE. The PINN's loss function consists of three error components,  for the prediction of the neural network as in the following  (i) Initial Condition, (ii) Boundary Condition, and (iii) PDE. Let $\hat{h}(x,t)$ be the output of neural network. The three components of the PINN's loss function are given below: 
\begin{itemize}
    \item Mean squared error on the Initial Condition
    \begin{align}
        \textrm{MSE}_I = \frac{1}{N_i}\sum_{k=1}^{N_i} \left(\hat{h}(x_{k}^{i},0) - h_{k}^{i}\right)^{2}\,, \quad x_{k}^{i} \in \Omega 
    \label{eq:3}
    \end{align}
    where $\hat{h}(x_{k}^{i},0)$ is the neural network output and $h_{k}^{i}$ is the given initial condition at $(x_{k}^{i},0)$. Here, the superscript, $(\bullet)^{i}$ stands for initial condition. 
     \item Mean squared error on the Boundary Condition
    \begin{align}
    \begin{gathered}
        \textrm{MSE}_B = \frac{1}{N_b}\sum_{k=1}^{N_b}\sum_{d=1}^{n_d}   \left(\hat{h}^{(d-1)}(x_{k}^{b},t_{k}^{b}) - \hat{h}^{(d-1)}(-x_{k}^{b},t_{k}^{b})\right)^{2}\;,
        \qquad (x_{k}^{b},t_{k}^{b}) \in \Gamma \times (0,T]
    \label{eq:4}
    \end{gathered}
    \end{align}  
    where $n_d$ is the highest order of derivative to which the periodicity is enforced on the boundary, $\Gamma$. Here, the superscript, $(\bullet)^{b}$ stands for boundary condition. 
    \item The Mean squared error due to Residual of the partial differential equation
    \begin{align}
    \begin{gathered}
        R := \hat{h}_t - F(\hat{h}, \hat{h}_x^{(1)},  \hat{h}_x^{(2)},...\hat{h}_x^{(m)}) \\
        \textrm{MSE}_R = \frac{1}{N_r}\sum_{k=1}^{N_r} \left(R(x_{k}^{r},t_{k}^{r})\right)^{2}\;,
        \qquad (x_{k}^{r},t_{k}^{r}) \in \Omega \times (0,T]
    \label{eq:5}
    \end{gathered}
    \end{align}
\end{itemize}
The superscript, $(\bullet)^{r}$ stands for residual of the PDE. ($x_k^{i}$) and (${x_k^{b},t_k^{b}}$), represent the set of points where the initial and boundary errors are computed. The residual/collocation error is computed at the collocation points (${x_k^{r},t_k^{r}}$). These points on the domain and the boundary are obtained using a latin hypercube sampling approach. 
Therefore, the total loss function of the neural network is given by adding all the aforementioned mean squared errors 
\begin{align}
    \textrm{MSE} = \textrm{MSE}_I + \textrm{MSE}_B + \textrm{MSE}_R
    \label{eq:6}
\end{align}

Once the PINN  is trained, the accuracy of the predicted solution is computed with respect to the true/exact solution at unknown  points (called testing points). Highly accurate solution of the initial boundary value problem obtained by the Chebyshev polynomial based numerical algorithm  \cite{Trefethen2014} and is considered as the exact solution. The relative total error ($\varepsilon_{total}$) of the PINN's prediction over the entire domain is obtained by normalizing the error with respect to the true solution as
\begin{align}
        \varepsilon_{total} = 
        \frac{\left[\frac{1}{N}\sum_{k=1}^{N} \left(\hat{h}(x_{k},t_{k}) - h(x_{k},t_{k})\right)^2 \right]^{1/2}}{\left[\frac{1}{N}\sum_{k=1}^{N} \left(h(x_{k},t_{k})\right)^2 \right]^{1/2}} 
        \label{eq:7}
\end{align}

The relative error ($\varepsilon$) of the PINN's prediction at each point  is obtained by normalizing the absolute error with respect to the true solution as
\begin{align}
        \varepsilon (x_{k},t_{k}) = 
        \frac{\abs{\hat{h}(x_{k},t_{k}) - h(x_{k},t_{k})}}{\left[\sum_{k=1}^{N} \left(h(x_{k},t_{k})\right)^2 \right]^{1/2}} 
        \label{eq:7a}
\end{align}

Where $h(x_{k},t_{k})$ is the true solution and $\hat{h}(x_{k},t_{k})$ is the neural network prediction for a set of testing points $\left\{(x_{k},t_{k})\right\}_{k=1}^{N}$, $ (x_{k},t_{k}) \in \Omega \times (0,T]$. For all comparisons between true and predicted solutions the relative total error `$\varepsilon_{total}$' and relative error `$\varepsilon$' is used.

\section{The proposed backward compatible sequential PINN method (BC-PINN)}\label{sec:BC-PINN}
In this section, we introduce an extension of the standard PINN technique that solves an initial-boundary value problem sequentially in time.
\subsection{bc-PINN}
\begin{figure}[h!]
    \centering
    \hspace{-5 mm}
    \includegraphics[scale=0.6]{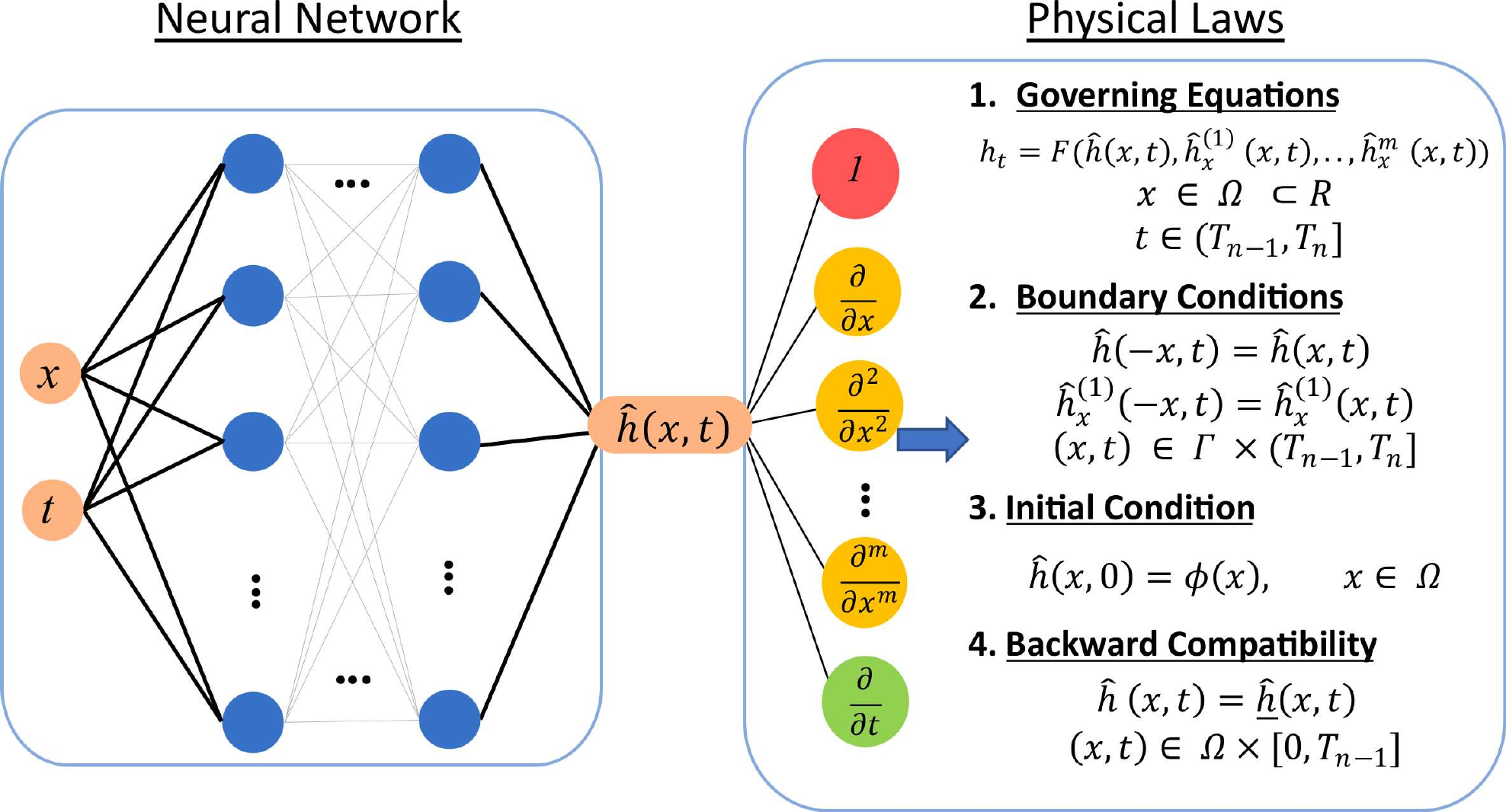}
    \caption{The schematics of the proposed backward compatible PINN (bc-PINN) approach for a time segment ($(T_{n-1},T_{n}]$). The neural network re-trains the PDE over $(T_{n-1},T_{n}]$ while satisfying the solution for all previous time segments.}     \label{fig:bc-PINN_schematics}
\end{figure}

\begin{figure}[h!]
    \centering
    \includegraphics[scale=0.55]{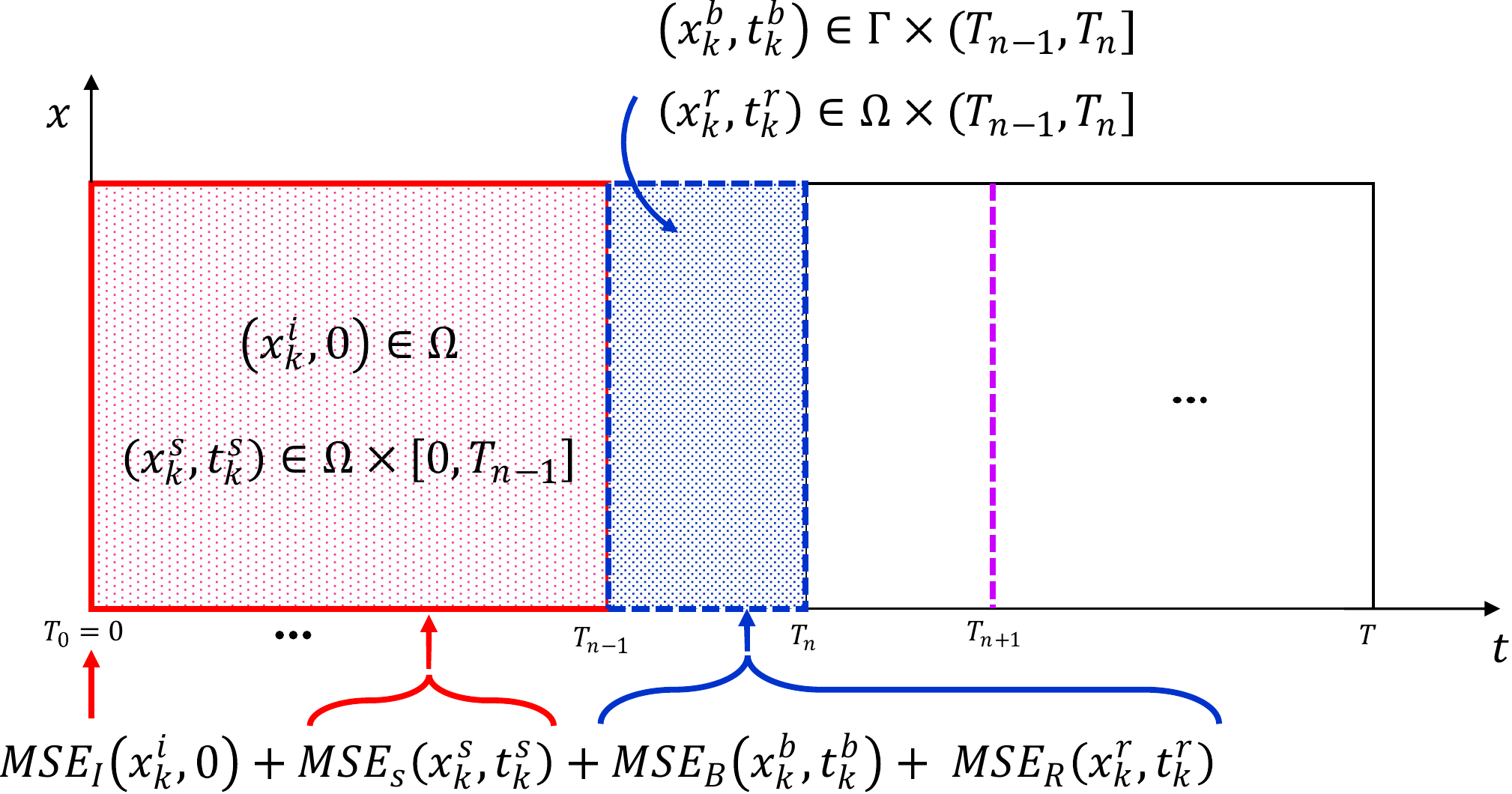}
    \caption{Illustration of the proposed backward compatibility scheme that satisfies the solutions obtained on all previous time segments  ($[0,T_{n-1}]$) while satisfying the PDE on the current time segment ($(T_{n-1},T_n$]).} 
    \label{fig:bc-PINN_sequential_scheme}
\end{figure}

 In the proposed method the PDE is solved progressively in time by re-training a single neural network over successive time segments.
The limitation of such retraining is that the network can predict only for the latest time segment and cannot predict for previous time segments for those it has been trained earlier. To overcome this limitation, the proposed model is designed to satisfy the solution of all the previous time segments while solving the PDE over a particular time segment. This scheme ensures backward compatibility of the solution by a single network. The proposed method is henceforth referred as backward compatible PINN  (bc-PINN).  The schematics of bc-PINN for a particular time segment is shown in figure~(\ref{fig:bc-PINN_schematics}) and the sequential scheme of proposed bc-PINN approach is shown in figure~(\ref{fig:bc-PINN_sequential_scheme}).
In bc-PINN the time domain $[0,T]$ is discretized into $n_{max}$ segments as
\begin{equation}
[T_0=0,T_1],\; [T_1,T_2],\cdots,\; [T_{n-1},T_{n}],\;\cdots, [T_{n_{max}-1},T_{n_{max}}=T]    \label{eq:8}
\end{equation}
where the $n^{\mathrm{th}}$ segment is denoted as $\Delta T_n = [T_{n-1}, T_{n}]$, $n=1,\cdots,n_{max}$. 

For the first time segment $\Delta T_1$ the solution of the PDE is sought through the standard PINN by minimizing the following loss function     
\begin{align}
\begin{gathered}
     \textrm{MSE}_{\Delta T_1} =  \textrm{MSE}_I(x_k^{i},0) + \textrm{MSE}_B(x_{k}^{b},t_{k}^{b}) + \textrm{MSE}_R(x_{k}^{r},t_{k}^{r})  \\ \\ x_{k}^{i} \in \Omega, \qquad
    (x_{k}^{b},t_{k}^{b}) \in \Gamma \times (0,T_{1}] \qquad 
    (x_{k}^{r},t_{k}^{r}) \in \Omega \times (0,T_{1}]
    \label{eq:9}
    \end{gathered}
\end{align}
Here, (${x_k^{i},t_k^{i}} $) represent the set of points where the error on initial condition is computed and (${x_k^{b},t_k^{b}}$) represent the set of points where the error on boundary condition is computed within the time segment $\Delta T_1=(0,T_{1}]$. 
For all of the subsequent time segments (i.e. ${\Delta T_n},\, n=2,\cdots,n_{max}$) we propose a novel loss function, which satisfies  the solution of all previous time segments. The solution of all previous time segments is enforced  by penalizing the departure from the already obtained solutions from the previous training, as given by 

\begin{equation}
    \begin{alignedat}{2}
        \textrm{MSE}_{\Delta T_n} =  \textrm{MSE}_I(x_k^{i},0) + \textrm{MSE}_B(x_{k}^{b},t_{k}^{b}) + \textrm{MSE}_R(x_{k}^{r},t_{k}^{r}) \\
        \; \; + \textrm{MSE}_S(x_{k}^{s},t_{k}^{s})  \,, 
        \quad n=2,\cdots,n_{max} \\ 
        x_{k}^{i} \in \Omega, \qquad
    (x_{k}^{b},t_{k}^{b}) \in \Gamma \times (T_{n-1},T_{n}] \qquad \\ 
    (x_{k}^{r},t_{k}^{r}) \in \Omega \times (T_{n-1},T_{n}], \qquad 
    (x_{k}^{s},t_{k}^{s}) \in \Omega \times [0,T_{n-1}]
    \end{alignedat}
    \label{eq:10}
\end{equation}
\noindent Here, ($x_k^{i},0$) represent the set of points where the error on initial condition is computed and  (${x_k^{b},t_k^{b}}$) represent the set of points where the error on boundary conditions is computed within the time segment $(T_{n-1},T_{n}]$.
The residual/collocation error as given in equation (\ref{eq:5}) is computed at the collocation points (${x_k^{r},t_k^{r}}$). 
We also minimize the departure from the already obtained solution that were stored at the grid points  $(x_{k}^{s},t_{k}^{s})$. The solution obtained (on $(0,T_{n}]$) at the $n^{\mathrm{th}}$ segment is stored for using it in the $(n+1)^{\mathrm{th}}$  segment.\\ 

\subsection{Details of the neural network of bc-PINN}\label{sec:NN}
We have used a standard (deep) neural network with two input neurons consisting of the spatial variable ($x$) and temporal variable ($t$). The output of the neural network ($\hat{h}(x,t)$) approximates the solution of the PDE ($h(x,t)$). To avoid model bias due to input features of different scales we have performed ``min-max'' normalization to scale the data uniformly. The neural network has 4 hidden layers consisting 200 neurons in each layer with a $tanh$ activation function. The neural network has more than 100,000 learning parameters which have been initialized using the ``xavier initialization'' technique. The optimization of the loss function and updating the learning parameters (weights and biases of the neural network) is performed using the ADAM and LBFGS optimizers. The learning rate for ADAM optimizer is considered as 0.001 with all other parameters as suggested in \cite{kingma2017adam}. Following standard PINN, after training the neural network using the ADAM optimizer we again train it using the L--BFGS optimizer until one of the following stopping criteria is met: (\RN{1}) Maximum iterations are equal to 50,000 (\RN{2}) Maximum number of function evaluations are equal to 50,000 (\RN{3})  Maximum number of line search steps (per iteration) equal to 50 (\RN{4}) The maximum number of variable metric corrections used to define the limited memory matrix are equal to 50 (\RN{5}) The iteration stops when $\frac{f^k - f^{k+1}}{max{(|f^k|,|f^{k+1}|,1)}} <= 2.22044604925e-16$, where $f$ is the neural network objective function and $k$ is the iteration number. 

\subsection{Details of the Computational Platform}\label{sec:Computer}
All the neural networks are trained on Nvidia Tesla P100 (3584 CUDA cores and 16GB of HBM2 vRAM) and Nvidia Volta V100 GPU (5120 CUDA cores, 640 Tensor cores  and 16GB of HBM2 vRAM). For inferencing and generating the true solutions via chebfun, we have used Dell precision 3630 workstation with Intel core i7-9700k 8 core (4.9 GHz Turbo) and 32 GB RAM. The software packages used for all the computations are Tensorflow 1.15 and MATLAB R2020a. \\

\subsection{The reference solution}\label{sec:exactSolution}
Accurate numerical solutions for the Allen Cahn and Cahn Hilliard equations are obtained using the chebfun package \cite{Trefethen2014}.
The  chebfun approach provides a polynomial interpolant for smooth functions in Chebyshev points.  
To solve time varying PDEs an exponential time differencing with Runge--Kutta time stepping scheme \cite{COX2002430} has been implemented in chebfun, which is used in the present work.  Henceforth, these solutions are considered as the true/exact solutions. We have taken 512 points for spatial discretization and 201 points for discretization in time scale. A fourth order Runge--Kutta time  integrator with time step $\Delta t = 10^{-5} $ is used.\\

The bc-PINN approach is applied to solve Allen Cahn equation and Cahn Hilliard equation in the next section to demonstrate its advantages for nonlinear and higher order PDEs in comparison to standard PINN method \cite{Raissi2019}.

\section{Allen Cahn Equation}\label{sec:AC}
\subsection{Allen Cahn equation and parameters}\label{sec:AC_equation}
The Allen Cahn equation is a semilinear partial differential equation which is well known for certain phase separation problems \cite{ALLEN19791085,bartels2015numerical,shen2010numerical}. For every $x \in \Omega$, ($\Omega$ is an open set of $\mathbb{R}^n$) the Allen Cahn equation is the $L^2$ gradient flow of the functional 
\begin{equation}
    I_{c_1}(h) = \frac{1}{2}\int_{\Omega} |\nabla h|^2 dx + \frac{1}{c_{1}^{2}}\int_{\Omega} F(h) dx
\end{equation}
For a phase separation problem, the parameter $h$ represents the concentration of the individual component and the parameter $c_1$ represents the interfacial thickness. The solution $h$ tends to minimize the energy functional $I_{c_1}(h)$ such that the it tends to the minimums of $F(h)$. The solution progressively  develops interfaces separating different phases. 
For a given initial condition, $h_0 \; \in \; L^2(\Omega)$ and $T \; > \; 0$ we seek a function $h: \Omega \times (0,T] \to \mathbb{R}$ which satisfies the equation below.

\begin{align}
    \begin{gathered}
        h_t - c_{1}^{2} \; \nabla^{2}h + f(h)= 0\, , \quad         t \, \in \, (0,T],  \; \; x \, \in \, \Omega \subset \mathbb{R}  \\
        f(h) = c_{2}(h^{3} - h) \\ 
        h(x,0) = x^2\cos(\pi x) \\ 
        h(x,t) = h(-x,t) \\ 
        h_{x}^{(1)}(x,t) = h_{x}^{(1)}(-x,t)
    \end{gathered}
    \label{eq:13}
\end{align}

The function $f$ is the derivative with respect to $h$ of a double well potential function $F$ where, $F \in C^1(\mathbb{R})$ is a non-negative function satisfying $F(\pm 1) = 0$. We have considered $ F = 5(h^2-1)^2/4$ and $f = 5(h^3 - h)$. Here the values of parameters for equation (\ref{eq:13}) are $c_1^2 = 0.0001$ and $c_2 = 5$. \\

\subsection{PINN for Allen Cahn equation}\label{sec:PINN_AC}

 \begin{figure}[!htb]
     \centering
     \subfloat[]{\label{fig:AC_sPINN_contour}\includegraphics[width=1\linewidth,valign=t]{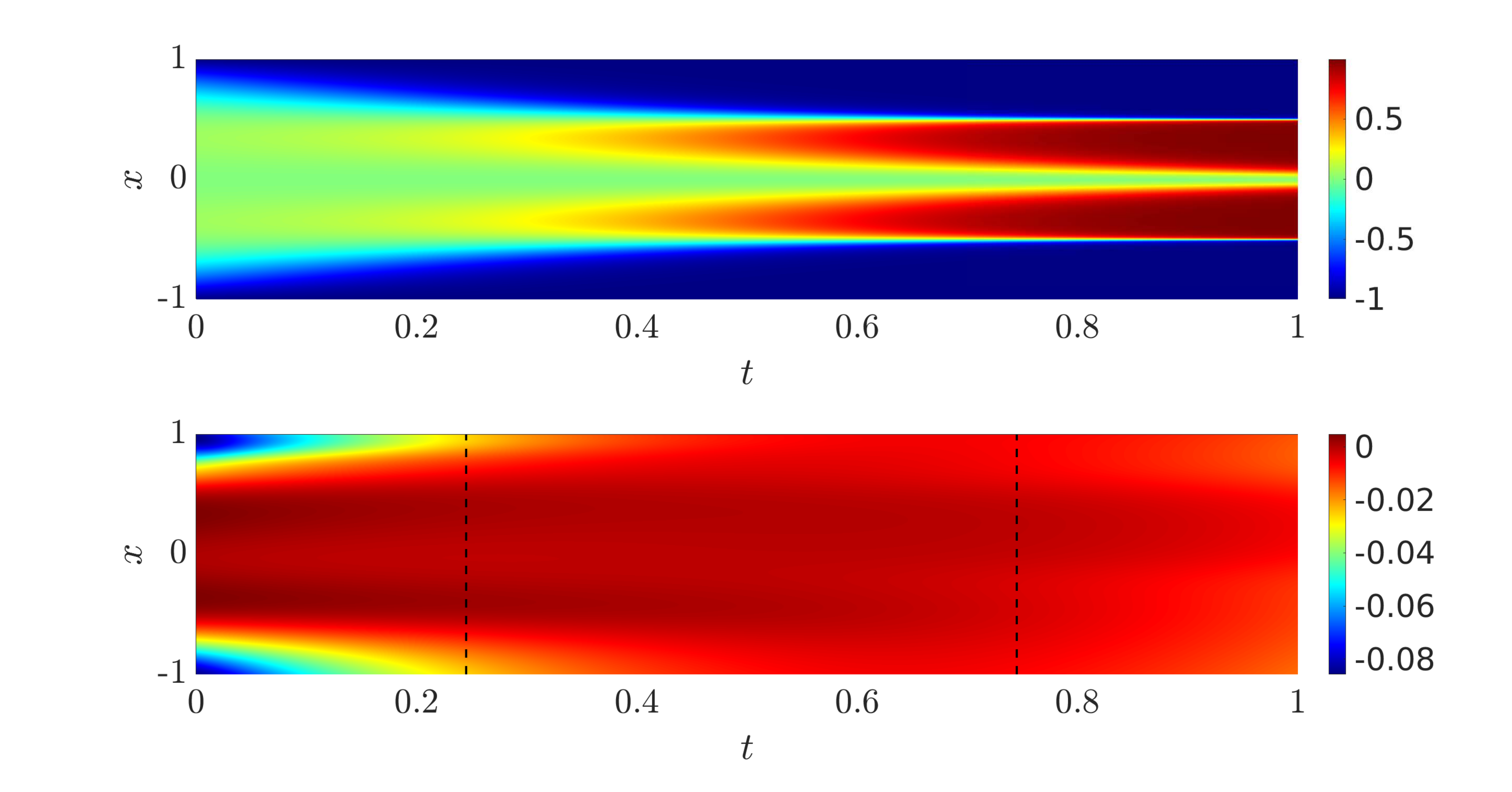}} \quad
     \subfloat[]{\label{fig:AC_sPINN_50_150}\includegraphics[width=0.8\linewidth,valign=t]{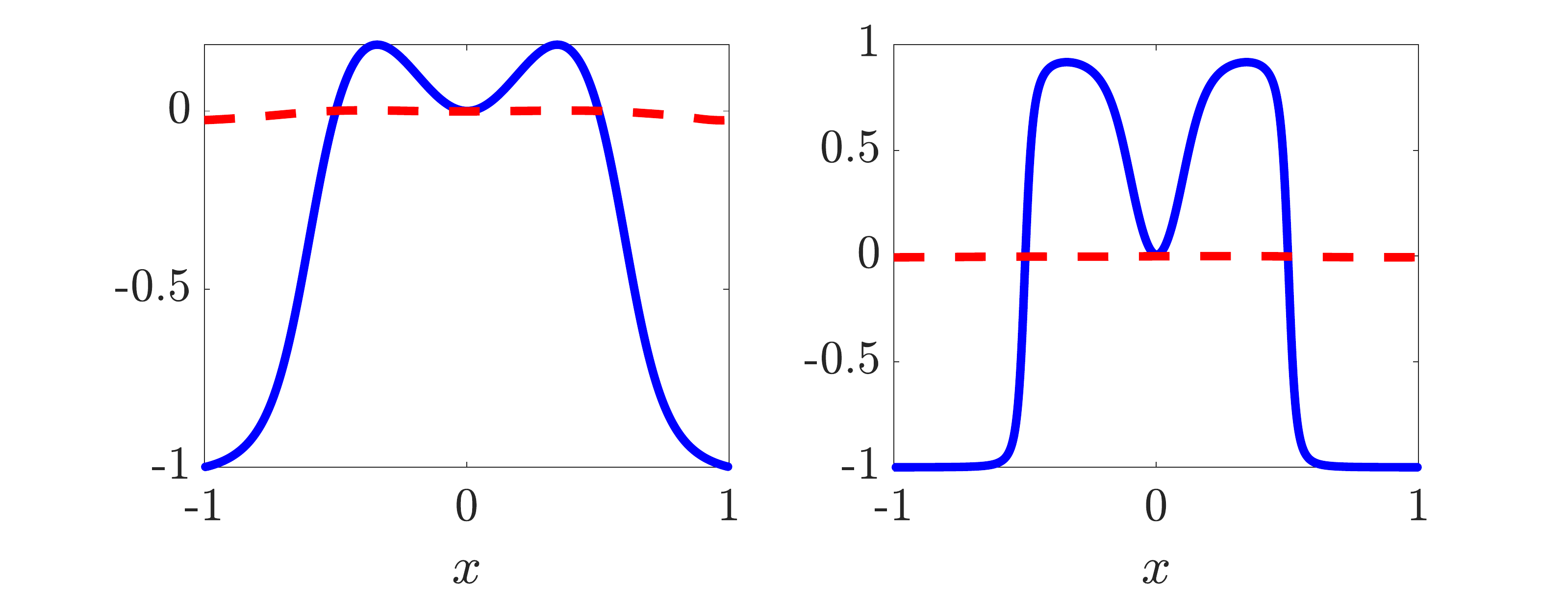}}
     \caption{(a): The exact solution (Top) and the PINN solution (Bottom) of the Allen Cahn equation for the entire spatio--temporal domain. (b): Time snapshots for the exact solution (\blueline) and the PINN solution (\redline) at $t=0.25$ and $t=0.75$.}
     \label{fig:AC_sPINN}
\end{figure}

\begin{figure}[!htb]
    \centering
    \includegraphics[scale=0.35]{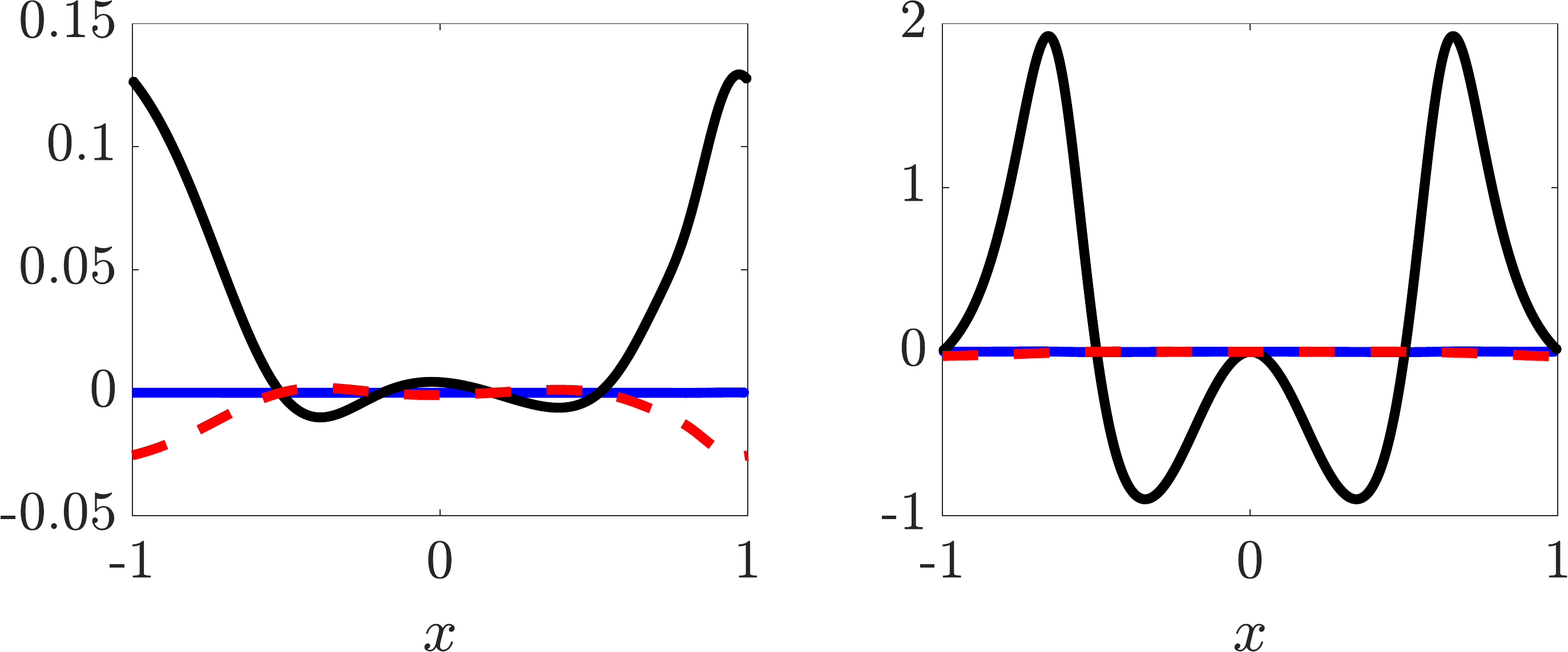}
    \caption{Individual terms of the Allen Cahn Equation obtained through the (Left): standard PINN method (Right): Chebfun method.  $h(x,t)$, ({\redline}), $0.0001\nabla^{2}h$ ({\blueline}) and $5(h^{3}-h)$ ({\blackline}) at t = 0.25.}
    \label{fig:AC_terms}
\end{figure}
At first we solve the Allen Cahn equation on $\Omega=[-1,1]$ and $(0,1]$ using standard PINN to demonstrate the challenge associated with non-linearity. For training the PINN, we have used 20,000 collocation points and trained for 100,000 ADAM iterations. The loss function for PINN is described in equation~(\ref{eq:3}),(\ref{eq:4}) and (\ref{eq:5}). The solution of standard PINN is quite erroneous as shown in figure~(\ref{fig:AC_sPINN}). In order to understand the reason for failure of the standard PINN, we analyze its prediction for the individual terms of the Allen Cahn equation. Figure~(\ref{fig:AC_terms}) shows the individual terms of the Allen Cahn equation obtained through the Chebfun method and the PINN. We observe that  the PINN fail to predict the non-linear term ($5(h^3-h)$) of the Allen Cahn equation.  
Therefore we have shown that the standard PINN \cite{Raissi2019} does not work for the Allen Cahn equation that consist of a strongly non-linear term. 

\subsection{bc-PINN for Allen Cahn equation}\label{sec:bcPINN_AC}
To overcome this limitation of PINN, we use the backward compatible PINN approach. In the bc-PINN  we propose the following loss function for any given time segment $\Delta T_{n}$ as given below.

\begin{itemize}
    \item Mean squared error on the initial Condition is the same as equation~(\ref{eq:3}). 
     \item Mean squared error on the boundary Condition
    \begin{align}
    \begin{gathered}
        \textrm{MSE}_B = \frac{1}{N_b}\sum_{k=1}^{N_b}\sum_{d=1}^{n_d} \left(\hat{h}^{(d-1)}(x_{k}^{b},t_{k}^{b}) - \hat{h}^{(d-1)}(-x_{k}^{b},t_{k}^{b})\right)^{2}\,,
        \qquad (x_{k}^{b},t_{k}^{b}) \in \Gamma \times (T_{n-1},T_{n}]
    \label{eq:15}
    \end{gathered}
    \end{align}  
    where $n_d$ is the order to which periodicity is enforced on the boundary $\Gamma$. Here, the superscript, $(\bullet)^{b}$ stands for boundary condition.
    \item The Mean squared error due to residual of the partial differential equation 
    \begin{align}
    \begin{gathered}
        R := \hat{h}_t - c_{1}^{2} \; \nabla^{2}\hat{h} + f(\hat{h})\\
        \textrm{MSE}_R =  \frac{1}{N_r}\sum_{k=1}^{N_r}   \left(R(x_{k}^{r},t_{k}^{r})\right)^{2}\,, 
        \qquad (x_{k}^{r},t_{k}^{r}) \in \Omega \times (T_{n-1},T_{n}]
    \label{eq:16}
    \end{gathered}
    \end{align}
    The superscript, $(\bullet)^{r}$ stands for residual of the PDE.
    \item Mean squared error for backward compatibility 
    \begin{align}
    \begin{gathered}
        \textrm{MSE}_S = \frac{1}{N_s}\sum_{k=1}^{N_s} \left(\hat{h}(x_{k}^{s},t_{k}^{s}) - \underline{\hat{h}}(x_{k}^{s},t_{k}^{s}) \right)^{2}\,,
        \qquad (x_{k}^{s},t_{k}^{s}) \in \Omega \times [0,T_{n-1}]
    \label{eq:17}
    \end{gathered}
    \end{align} 
    where, ${\hat{h}}(x,t)$ is the neural network prediction and $\underline{\hat{h}}(x,t)$ is the known solution through the neural network from the previous time steps $\Omega \times [0,T_{n-1}]$. The superscript, $(\bullet)^{s}$ stands for the backward compatible solution. 
    \item The total mean squared error or loss is given as
    \begin{align}
        \textrm{MSE}_{\Delta T_{n}} = \textrm{MSE}_I + \textrm{MSE}_B + \textrm{MSE}_R + \textrm{MSE}_S
    \label{eq:18}
    \end{align}
\end{itemize}

 The hyper-parameters associated with training the bc-PINN are number of ADAM iterations ($N_{iter}$), time steps per segment and number of collocation points ($N_r$) per segment. Table \ref{tab:AC_training_data} shows the values of all the aforementioned hyper-parameters for training the bc-PINN. \\
\begin{table}[!htb]
\centering
\begin{tabular}{|c|c|c|}
\hline
\textbf{Variable} & \textbf{Description} & \textbf{Number} \\ \hline
$N_i$ & Initial Datapoints & 512 \\ \hline
$N_b$ & Boundary Datapoints & 40/segment \\ \hline
$N_r$ & Collocation points & 20000/segment \\ \hline
$N_{iter}$ & Number of ADAM iterations & 10000/segment \\ \hline
\end{tabular}
\caption{\small{Description of Training Data for Allen Cahn Equation. The segment considered here consists of 40 time steps. }}
\label{tab:AC_training_data}
\end{table}

\begin{figure}[!htb]
     \centering
     \subfloat[Standard PINN]{\label{fig:AC_sPINN_t_50}\includegraphics[width=0.3\linewidth]{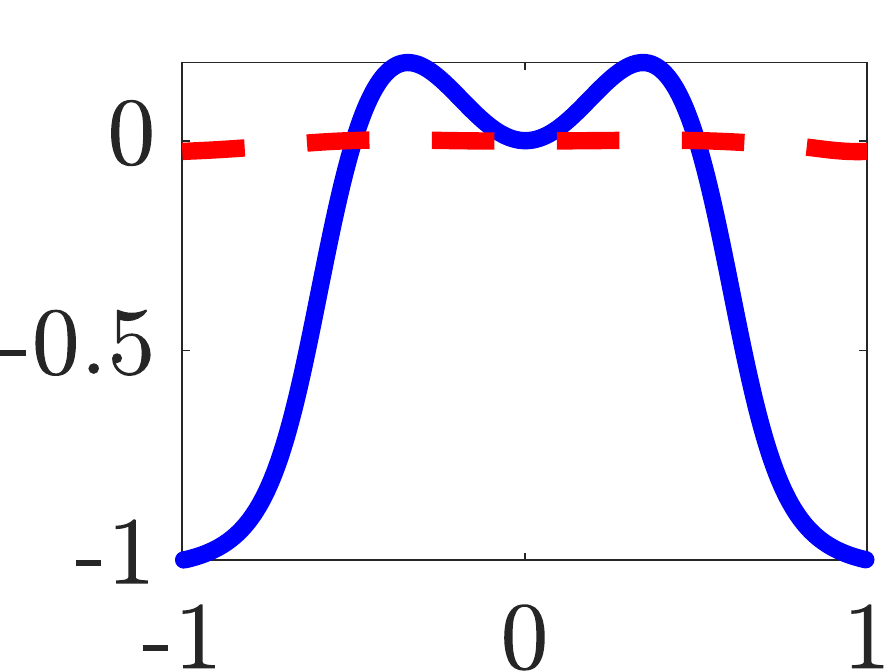}}\qquad
     \subfloat[bc-PINN]{\label{fig:AC_bcPINN_t_50}\includegraphics[width=0.3\linewidth]{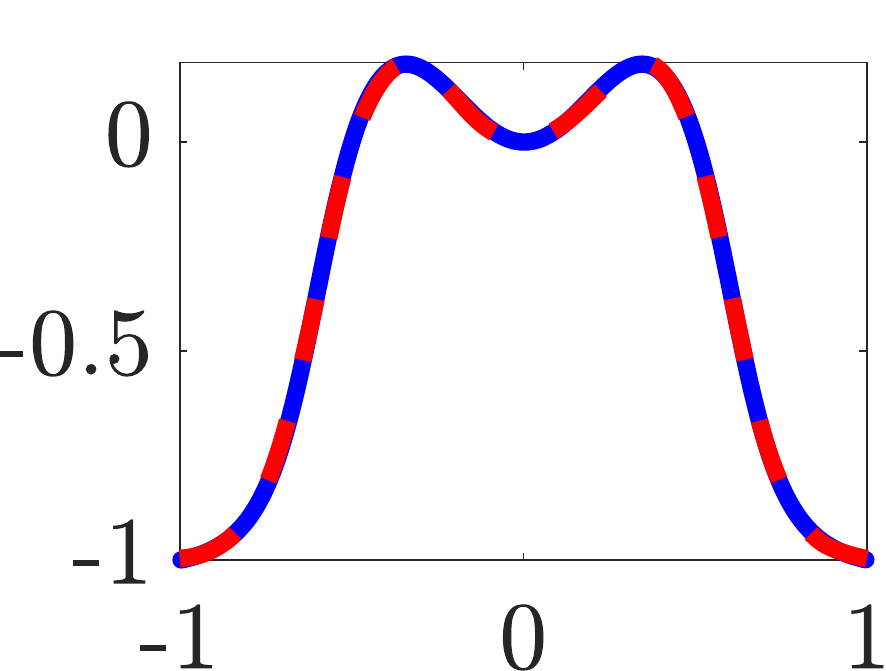}}
     \vspace{10pt}
     \caption{Exact (\blueline) and Predicted (\redline) solution at time t = 0.25}
     \label{fig:AC_bc-PINN-solutions}
\end{figure}

The exact and predicted solution at time $t=0.25$ obtained by the standard PINN and bc-PINN are shown in figure~(\ref{fig:AC_bc-PINN-solutions}). While the standard PINN fails, the proposed bc-PINN predicts the solution quite accurately. The relative total errors ($\varepsilon_{total}$) for both the approaches are shown in table~\ref{tab:AC_error}. 

\begin{table}[h!]
\centering
\begin{tabular}{|c|c|} 
\hline
\textbf{Method} & \textbf{Error($\mathbf{\varepsilon_{total}}$)} \\ \hline
Standard PINN & 0.9919  \\ \hline
bc-PINN & 0.0701 \\ \hline
\end{tabular}
\caption{\small{Relative total errors (equation~(\ref{eq:7})) over the entire domain with respect to Chebfun solution for different methods.}}
\label{tab:AC_error}
\end{table}
  
\begin{figure}[h!]
     \centering
     \subfloat[]{\label{fig:AC_wolog}\includegraphics[width=1\linewidth,valign=t]{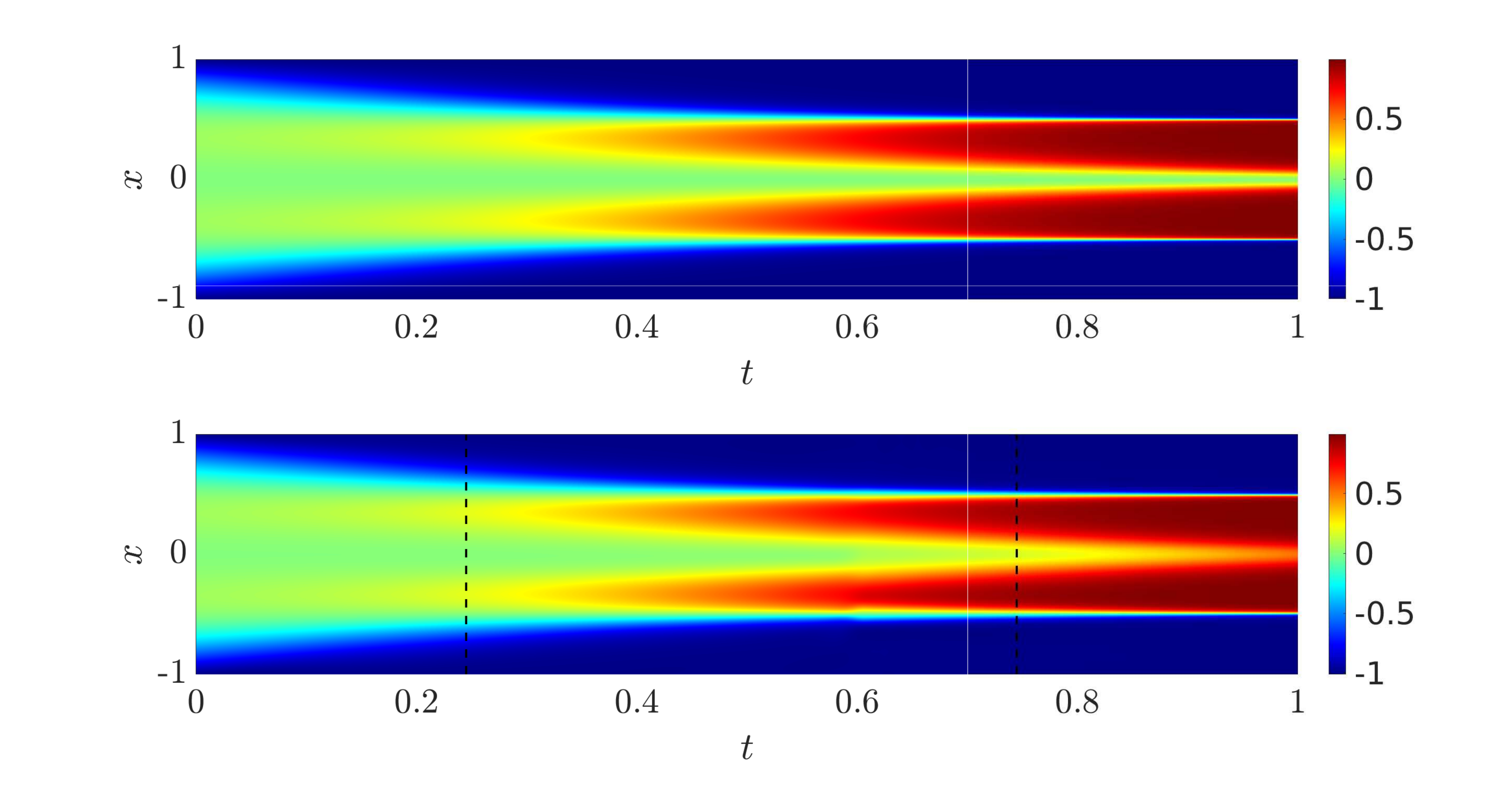}} \quad
     \subfloat[]{\label{fig:AC_50_150}\includegraphics[width=0.7\linewidth,valign=t]{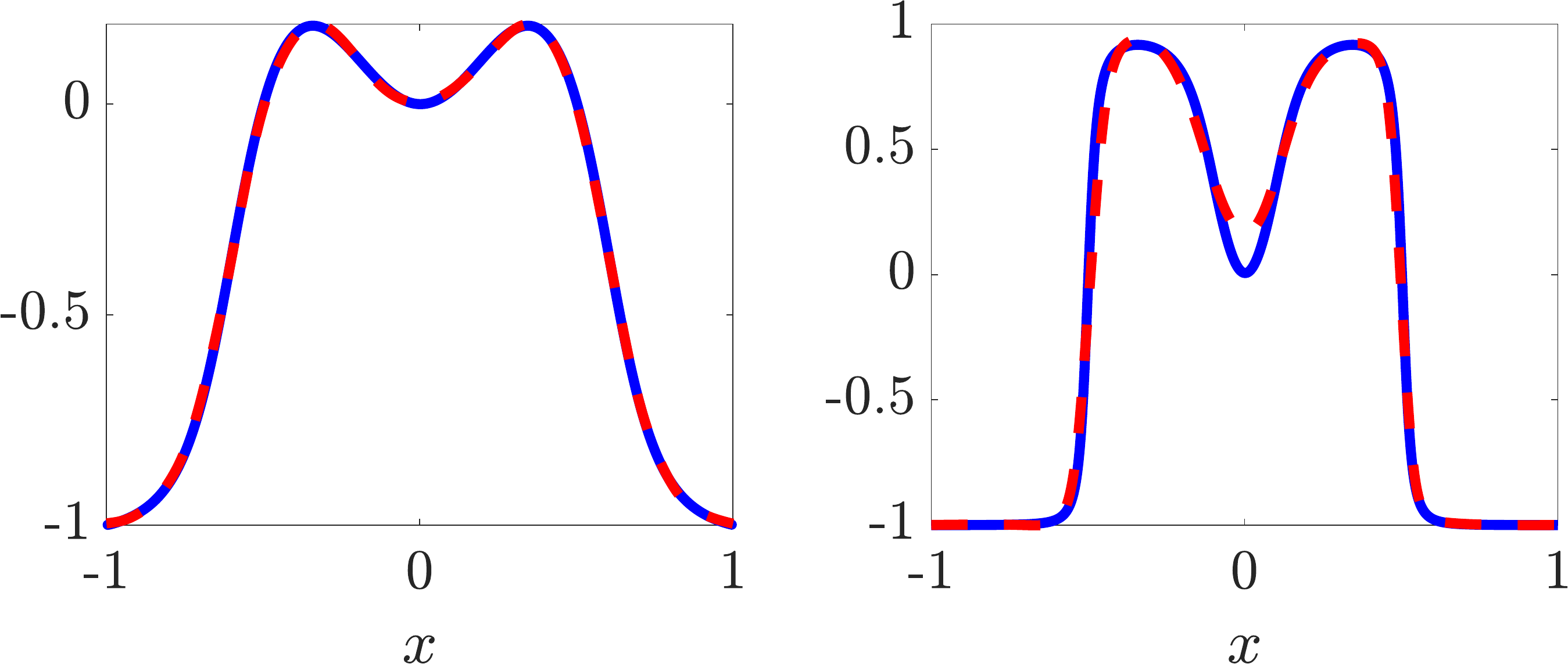}}
     \caption{(a): Exact (Top) and bc-PINN (Bottom) solutions of the Allen Cahn equation for the entire spatio--temporal domain. (b): The exact (\blueline) and the  bc-PINN (\redline) solutions at time $t=0.25$ and $t=0.75$.}
     \label{fig:AC_bc-PINN}
\end{figure}

\begin{figure}[htbp]
     \centering
     \subfloat[Predicted solution (top) and  relative error (bottom) obtained using PINN]{\includegraphics[width=1\linewidth]{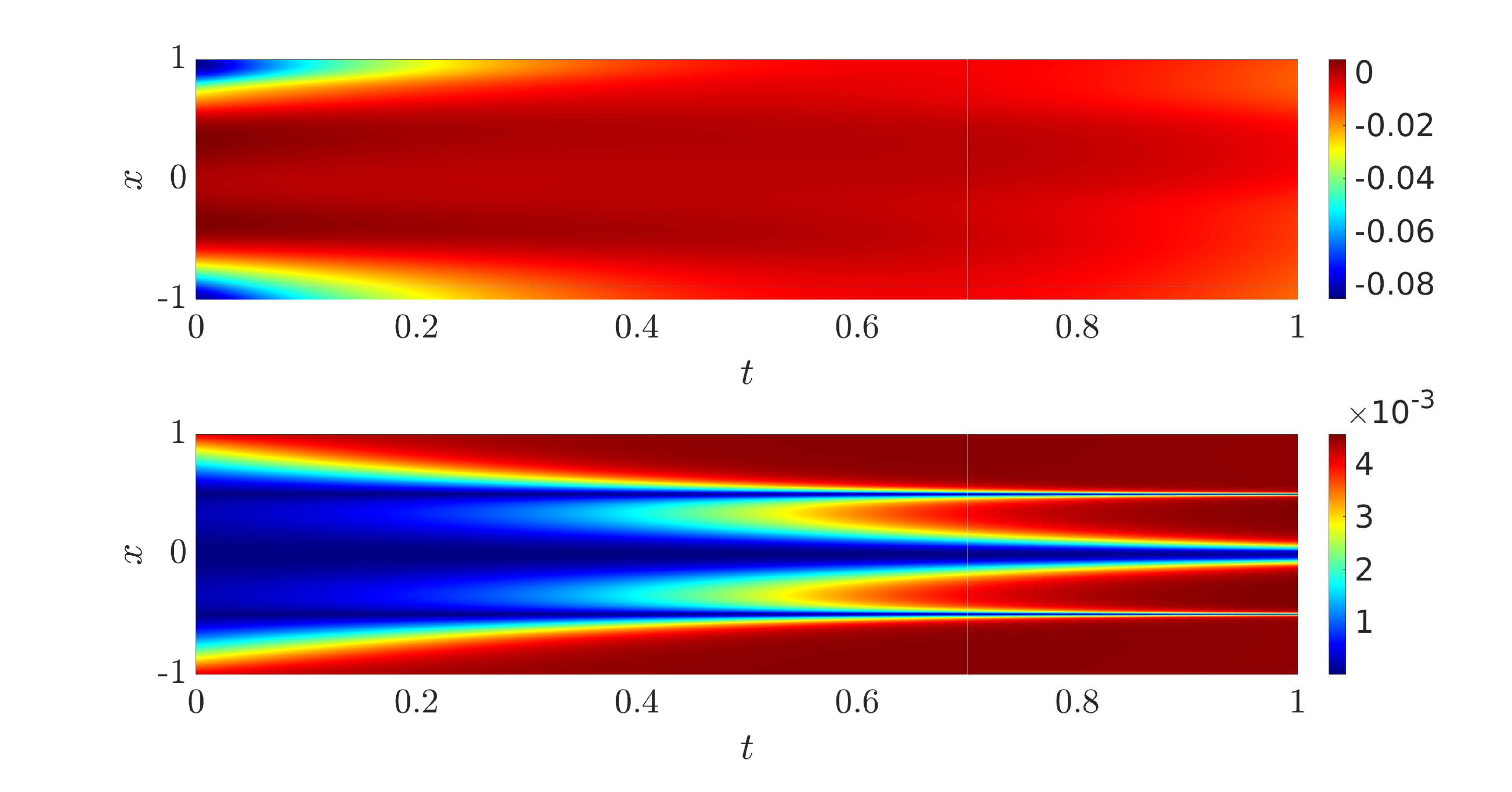}} \quad
     \subfloat[Predicted solution (top) and  relative error (bottom) obtained using bc-PINN]{\includegraphics[width=1\linewidth]{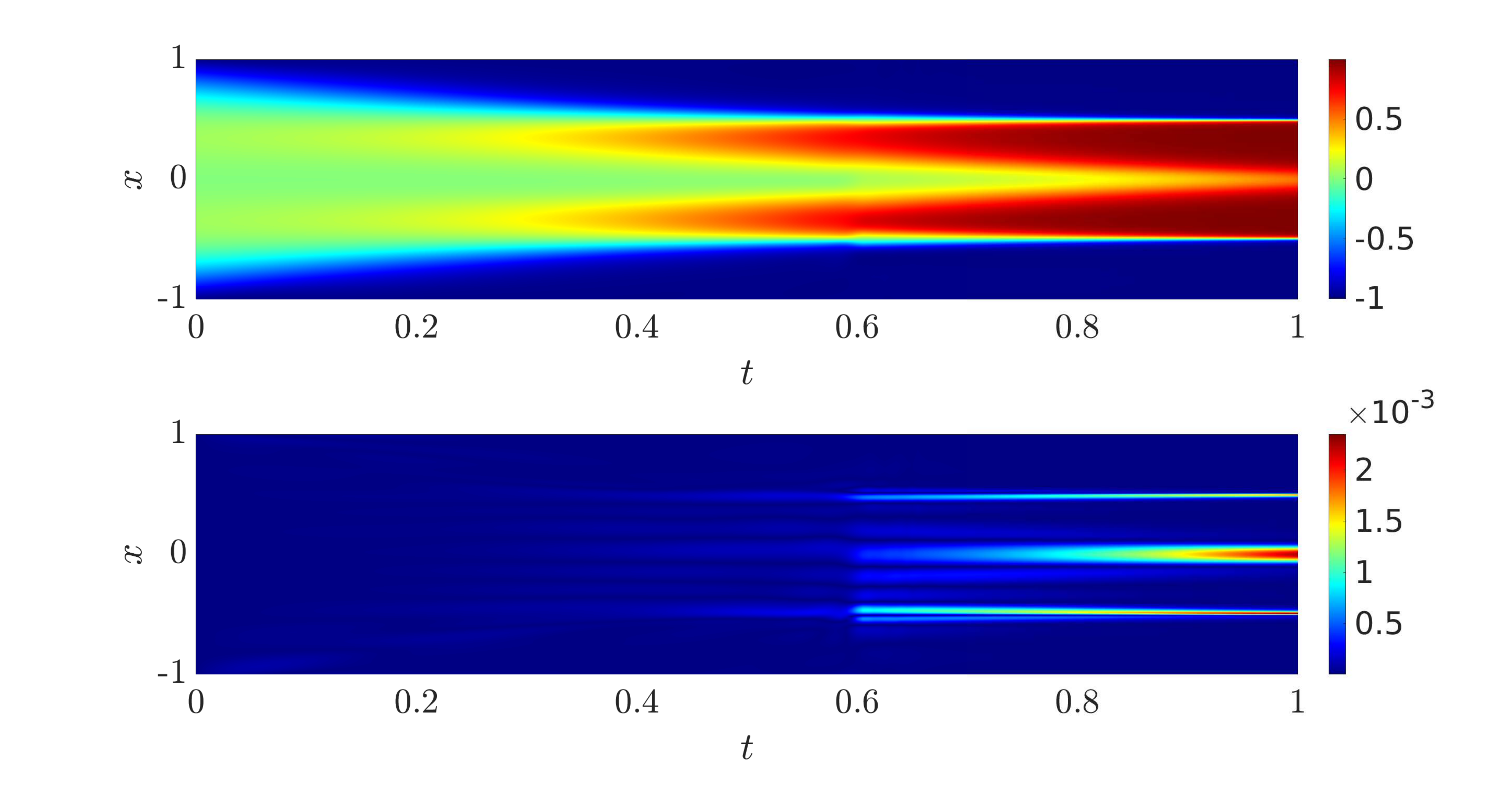}}
    \caption{Solution and relative error for Allen Cahn equation through bc-PINN.}
    \label{fig:AC_errplot}
\end{figure}

\begin{figure}[htbp]
     \centering
     \subfloat[Predicted solution (top) and  relative error (bottom) obtained using bc-PINN for $c_{1}^{2}$ = 0.00001]{\includegraphics[width=1\linewidth]{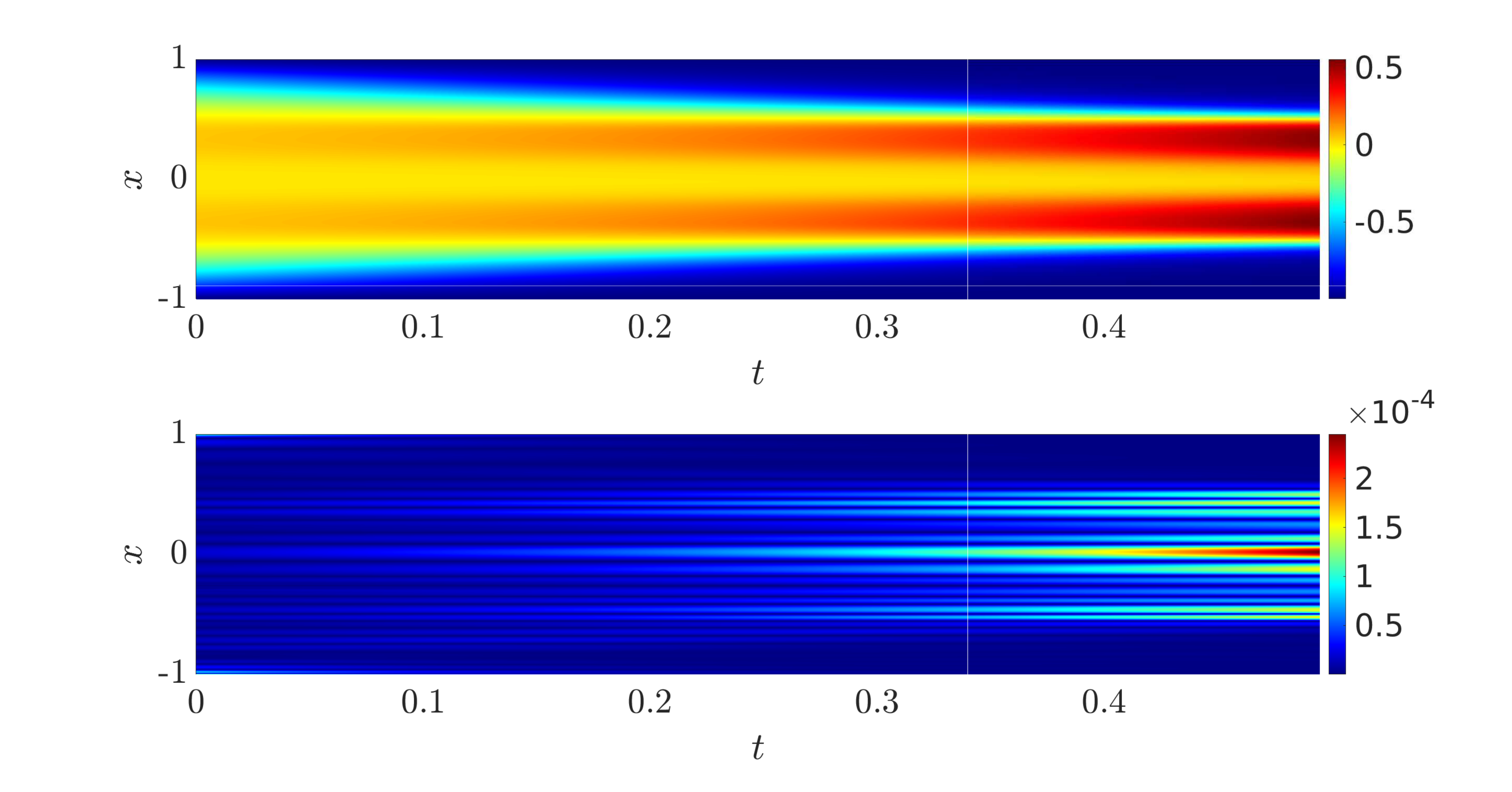}} \quad
     \subfloat[Predicted solution (top) and  relative error (bottom) obtained using bc-PINN for $c_{1}^{2}$ = 0.00005]{\includegraphics[width=1\linewidth]{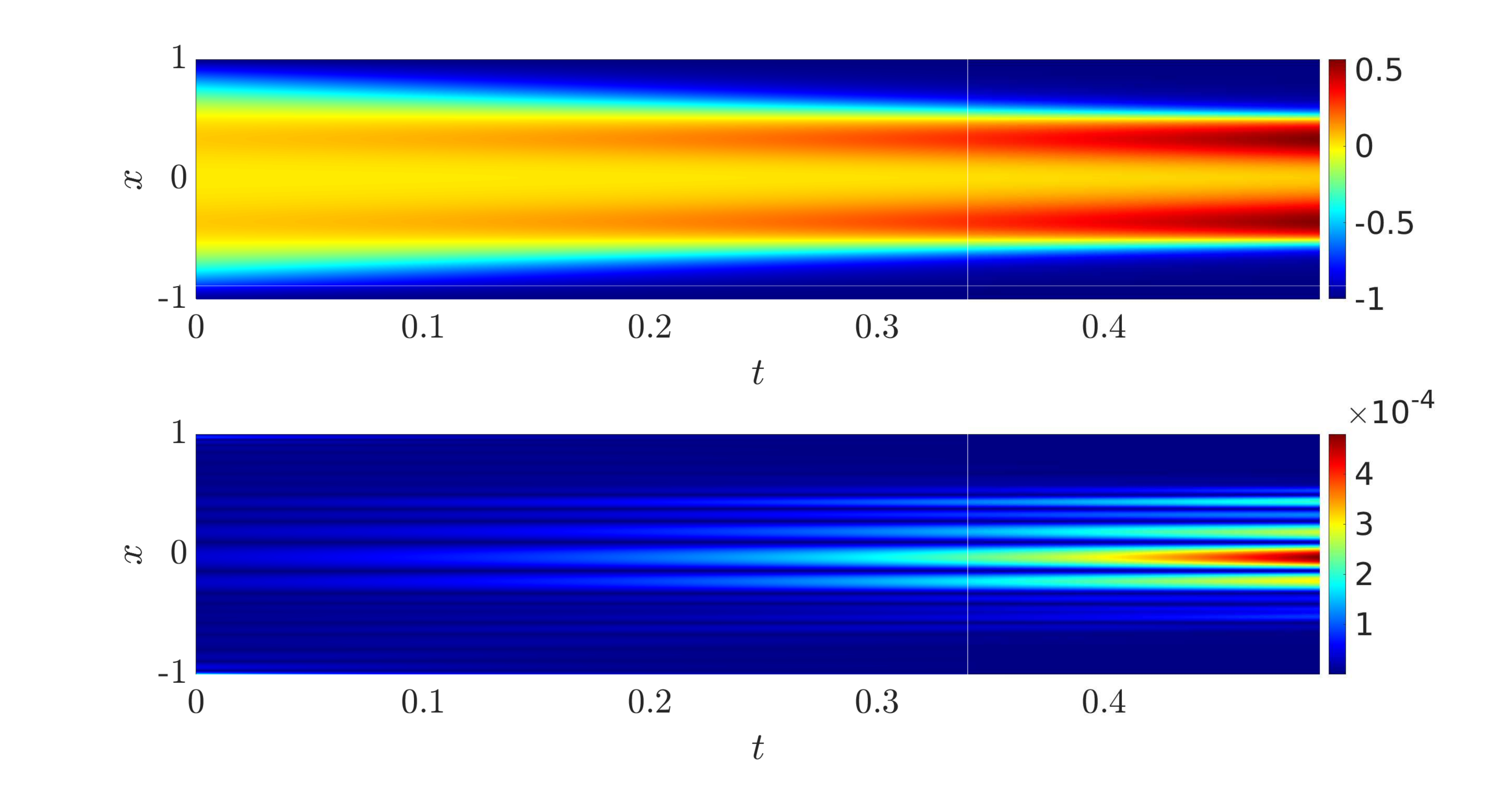}}
     \caption{Solutions and relative errors of the Allen Cahn equation for different  $c_{1}^{2}$ (of equation~(\ref{eq:13})) obtained by the bc-PINN method.} 
     \label{fig:AC_diff_C1}
\end{figure}

The comparison between the predicted solution using bc-PINN and the chebfun solution is shown in figure~(\ref{fig:AC_bc-PINN}). This shows that the bc-PINN  can accurately predict the solution for the entire domain. The solutions and errors by the PINN and bc-PINN are compared in figure~(\ref{fig:AC_errplot}), showing much higher accuracy by the bc-PINN. The error plots confirms high accuracy of bc-PINN. The error increases with time very slowly. This is due to two reasons: (\RN{1}) the solution becomes progressively phase-separated (between zero and one) yielding greater curvatures and sharp phase-boundaries that are difficult to capture, and (\RN{2}) due to the sequential nature of the bc-PINN approach the error accumulates with time progression, which is similar to the time-integrators. To illustrate the high accuracy of the  bc-PINN approach, solutions and errors for different values of the interfacial thickness ($c_1$) is plotted in figure~(\ref{fig:AC_diff_C1}). 
As we decrease the parameter $c_1$, it can be seen that the error in the prediction decreases. The parameter $c_1$ controls the effect of the double derivative of the solution ($\nabla^{2}\, h$). Therefore, as we decrease $c_1$ the error due to the approximation in derivative reduces and thus the accuracy of the bc-PINN solution increases. In ~\ref{sec:bc-PINN_logres} a new loss function including a logarithmic residual for the Allen Cahn equation is discussed. This new logarithmic residual bc-PINN approach and its results are presented in comparison with the simple bc-PINN approach without a logarithmic residual.

\clearpage
\section{Cahn Hilliard Equation}\label{sec:CH}
\subsection{Cahn Hilliard equation and parameters} \label{sec:CH_equation}
The Cahn Hilliard equation is considered here to demonstrate the advantages of the proposed bc-PINN method for nonlinear and higher order PDEs. 
The Cahn Hilliard equation \cite{cahn1958free,Alain_CH,kim2016basic} plays an essential role in the field of material science for describing the qualitative features in a phase separation process for two phase systems (assuming isotropy and constant temperature) \footnote{The process of phase separation can be observed when a binary alloy is cooled down adequately. This leads to a state of total nucleation which is mainly referred to as spinodal decomposition. In the subsequent stage coarsening occurs in the nucleated microstructure at a much slower rate. This whole phase separation phenomena affects the mechanical properties (eg. strength, hardness and fracture toughness) of the material.}. Since the entire process is governed by the Cahn Hilliard equation it is essential to understand the physical significance of each individual variable.
\begin{align}
  h_{t}(x,t) - \nabla^{2}(-\alpha \kappa \nabla^{2}h(x,t) + \kappa f(h(x,t))) = 0\, , \quad  x\in \Omega, \; \;  t \in (0,T]
\label{eq:20}
\end{align}
Here $\Omega$ is an open set in $\mathbb{R}$. 
The order parameter $h$ in equation~(\ref{eq:20}), refers to the rescaled density or concentration of one of the material components in the system and it takes values between (-1 and 1 which corresponds to their pure states). The density of second component is $1-h$, and this ensures that the total density over the simulation domain is a conserved quantity. In addition, the function $f(h)$, is the derivative of a double well potential function. The two phases of the system correspond to the two wells of the double well potential function. The parameter $\kappa$ is the mobility parameter and  the  parameter $\alpha$ is related to the surface tension at the interface. Here the values of parameters are taken as $\alpha = 0.02$ and $\kappa = 1$.


 \begin{figure}[!htb]
     \centering
     \subfloat[]{\label{fig:CH_sPINN_contour}\includegraphics[width=1.1\linewidth,valign=t]{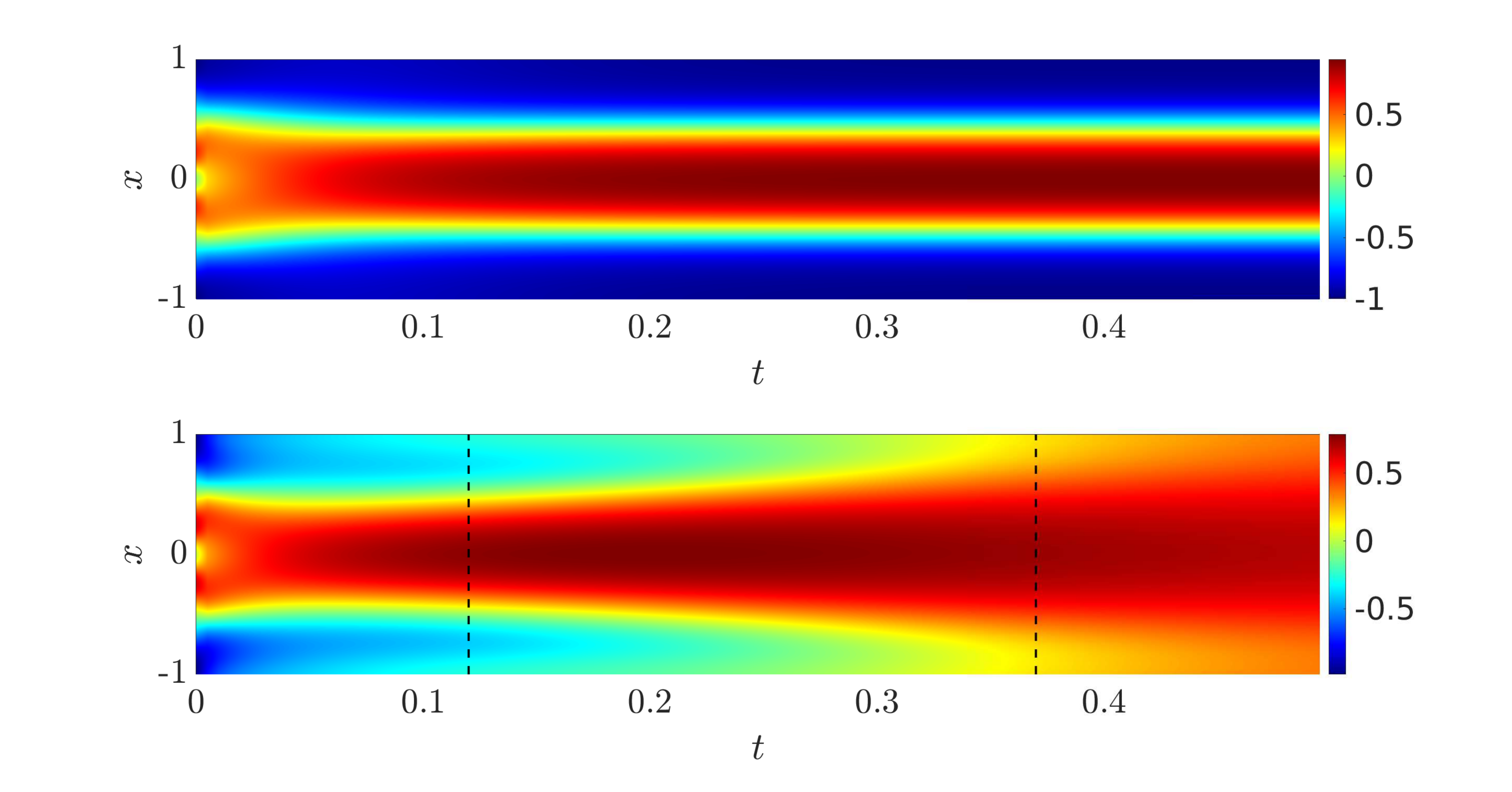}} \quad
     \subfloat[]{\label{fig:CH_sPINN_25_75}\includegraphics[width=0.85\linewidth,valign=t]{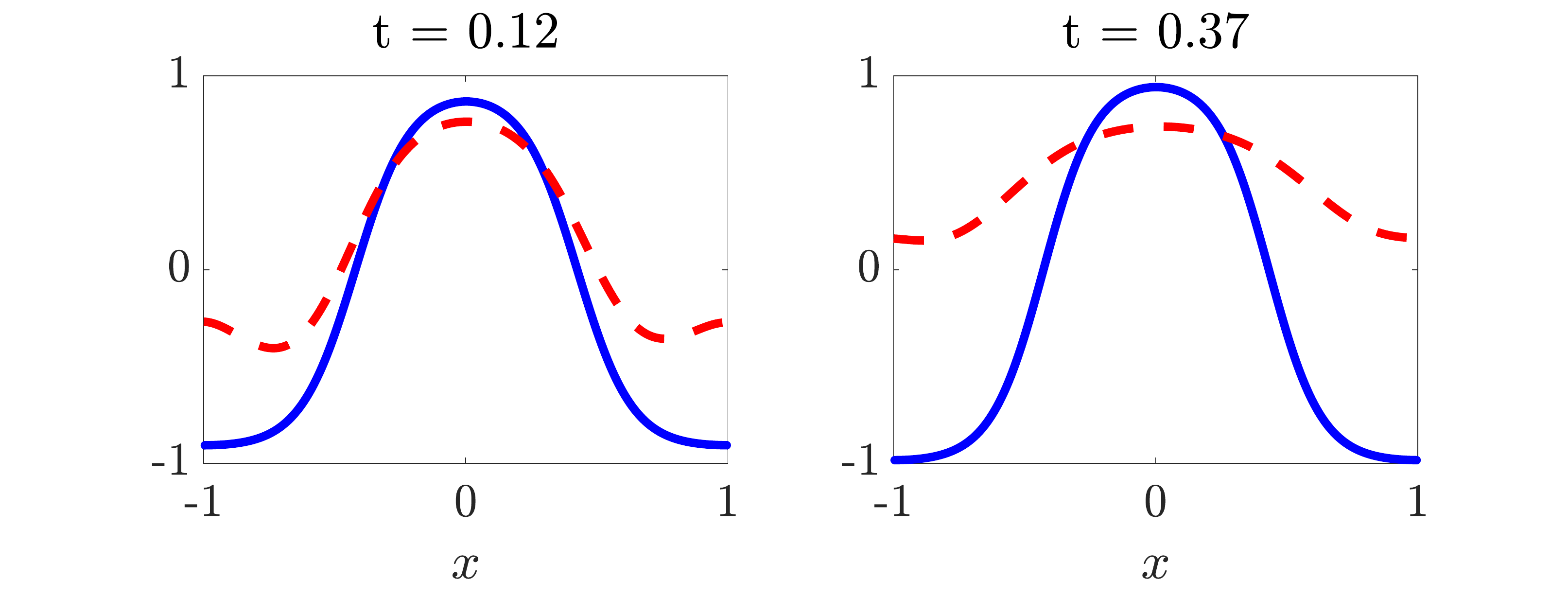}}
     \caption{(a): The exact solution (Top) and the PINN solution (Bottom) of the Cahn Hilliard equation for the entire spatio--temporal domain. (b): Time snapshots for the exact solution (\blueline) and the PINN solution (\redline) at $t=0.12$ and $t=0.37$.}
     \label{fig:CH_sPINN}
\end{figure}

\subsection{PINN for Cahn Hilliard equation}\label{sec:PINN_CH}
Initially we solve the Cahn Hilliard equation using standard PINN to demonstrate the challenge for this equation. The domain Cahn Hilliard equation on $\Omega = [-1,1]$ and periodic boundary conditions. The time domain considered is $(0,T] = (0,1]$. For training the standard PINN, we have used 20,000 collocation points and trained for 100,000 ADAM iterations. The loss function for PINN is described in equation~(\ref{eq:3}),~(\ref{eq:4}) and~(\ref{eq:5}). The solution predicted after training is shown in figure~(\ref{fig:CH_sPINN}) and it can be observed that there is significant mismatch between the PINN and the exact solution. \\

The two possible reasons for the inaccurate solution are strong non-linearity and  the high order derivative terms (fourth order). In PINN the derivatives are approximated using automatic differentiation.
It has been shown that as the order of the derivative increases the complexity in automatic differentiation increases and it becomes computationally expensive \cite{baydin2018automatic}. In order to overcome the difficulty in approximating the higher order derivative via automatic differentiation, we adopt the phase space representation in the proposed bc-PINN. 

\subsection{bc-PINN for Cahn Hilliard equation}\label{sec:bcPINN_CH}
In this section, we introduce the bc-PINN approach with a phase space representation for solving the Cahn Hilliard Equation. The phase space representation is widely used to represent a high order PDE into  coupled multiple lower order PDEs. We show that by adopting phase space representation the accuracy of the  bc-PINN method improves significantly.  
The phase space representation of the Cahn Hilliard equation (a fourth order PDE,  equation~(\ref{eq:20}))   yields two coupled second order PDEs. 
\begin{align}
    \begin{gathered}
        h_{t} = \nabla^{2}(-\alpha\, \kappa\,\mu + \kappa f(h)), \quad \mu = \nabla^{2}h \qquad
        t \, \in \, (0,T],  \; \; x \, \in \, \Omega \subset \mathbb{R}  \\
        h(x,0) = cos(\pi x) - \exp(-4 (\pi x)^2 )\\
        h(-x,t) = h(x,t), \quad (x,t) \in \Gamma \times (0,T]  \\
        h_x^{(1)}(-x,t) = h_x^{(1)}(x,t)\, , \quad (x,t) \in \Gamma \times (0,T]\\
        \mu(-x,t) = \mu(x,t), \quad (x,t) \in \Gamma \times (0,T]  \\
        \mu_x^{(1)}(-x,t) = \mu_x^{(1)}(x,t)\, , \quad (x,t) \in \Gamma \times (0,T] 
    \end{gathered}
    \label{eq:21}
\end{align}  
Where, $f(h) = h(h^{2} - 1)$. 
Therefore there are two outputs of the neural network $\hat{h}(x,t)$ and $\hat{\mu}(x,t)$ in the present method. The input features are the spatio--temporal variables $(x,t)$. The modified loss function for the coupled phase space system includes an error on initial condition, error on the boundary conditions and error on the residual. In addition it will have the error for the backward compatibility.  
Therefore, the total loss function (equation~(\ref{eq:26})) for any time segment $\Delta T_{n}$ is sum of all the aforementioned errors given in   equation~(\ref{eq:22}--\ref{eq:25}).

\begin{itemize}
    \item Mean squared error on the initial condition for $h(x,t)$ and $\mu(x,t)$
    \begin{align}
        \textrm{MSE}_I = \frac{1}{N_i}\left\{\sum_{k=1}^{N_i} \left(\hat{h}(x_{k}^{i},0) - h_{k}^{i}\right)^{2} + \sum_{k=1}^{N_i} \left(\hat{\mu}(x_{k}^{i},0) - \mu_{k}^{i}\right)^{2} \right\} \,, \quad x_{k}^{i} \in \Omega 
    \label{eq:22}
    \end{align}
    Here, the neural network output is $\hat{h}(x_{k}^{i},0),\hat{\mu}(x_{k}^{i},0)$ and the given initial condition is $h_{k}^{i},\mu_{k}^{i}$ at $(x_{k}^{i},0)$. The superscript, $(\bullet)^{i}$ stands for initial condition. 
     \item Mean squared error on the boundary Condition
    \begin{align}
    \begin{gathered}
        \textrm{MSE}_{B_{h}} = \frac{1}{N_b}\left\{\sum_{k=1}^{N_b}\sum_{d=1}^{n_d} \left(\hat{h}^{(d-1)}(x_{k}^{b},t_{k}^{b}) - \hat{h}^{(d-1)}(-x_{k}^{b},t_{k}^{b})\right)^{2} \right\}\\
        \textrm{MSE}_{B_{\mu}} = \frac{1}{N_b}\left\{\sum_{k=1}^{N_b}\sum_{d=1}^{n_d} \left(\hat{\mu}^{(d-1)}(x_{k}^{b},t_{k}^{b}) - \hat{\mu}^{(d-1)}(-x_{k}^{b},t_{k}^{b})\right)^{2} \right\}\\
        \textrm{MSE}_B = \textrm{MSE}_{B_{h}} + \textrm{MSE}_{B_{\mu}} \,,
         \qquad (x_{k}^{b},t_{k}^{b}) \in \Gamma \times (T_{n-1},T_{n}]
    \label{eq:23}
    \end{gathered}
    \end{align}  
    where $n_d$ is the order to which periodicity is enforced. The superscript, $(\bullet)^{b}$ stands for boundary condition.
    \item Mean squared error on the Residual of the partial differential equation
    \begin{align}
    \begin{gathered}
        R_1 := \hat{h}_t - \nabla^{2}\hat{h}(-\mu + f(\hat{h})) \\
        R_2 := \hat{\mu} - \nabla^{2}\hat{h} \\
        \mathrm{MSE}_R = \frac{1}{N_r}\left\{\sum_{k=1}^{N_r} R_{1}\left(x_{k}^{r},t_{k}^{r}\right)^{2} + \sum_{k=1}^{N_r} R_{2}\left(x_{k}^{r},t_{k}^{r}\right)^{2} \right\}\,, 
        \qquad (x_{k}^{r},t_{k}^{r}) \in \Omega \times (T_{n-1},T_{n}]
    \label{eq:24}
    \end{gathered}
    \end{align}
    The superscript, $(\bullet)^{r}$ stands for residual of the PDE.
    \item Mean squared error for backward compatibility 
    \begin{align}
    \begin{gathered}
        \mathrm{MSE}_S = \frac{1}{N_s}\sum_{k=1}^{N_s} \left(\hat{h}(x_{k}^{s},t_{k}^{s}) - \underline{\Tilde{h}}(x_{k}^{s},t_{k}^{s}) \right)^{2}\,,
        \qquad (x_{k}^{s},t_{k}^{s}) \in \Omega \times [0,T_{n-1}]
    \label{eq:25}
    \end{gathered}
    \end{align} 
    where, $\hat{h}(x,t)$ is the neural network prediction and $\underline{\hat{h}}(x,t)$ is the known solution through the neural network from the previous time steps $\Omega \times [0,T_{n-1}]$. The superscript, $(\bullet)^{s}$ stands for the backward compatible solution. 
    \item The total mean squared error or loss is given as
    \begin{align}
        \mathrm{MSE}_{\Delta T_{n}} = \mathrm{MSE}_I + \mathrm{MSE}_B + \mathrm{MSE}_R + \mathrm{MSE}_S
    \label{eq:26}
    \end{align}
\end{itemize} 

The boundary loss (equation~(\ref{eq:23}) is applied for $n_d = 1$ on  $\hat{h}$ and $\hat{\mu}$, to represent  periodic boundary conditions. Equation~(\ref{eq:24}) describes two components of residual for two PDEs in the phase space form of the Cahn Hilliard equation (equation~(\ref{eq:21})). 
\begin{table}[h!]
\centering
\begin{tabular}{|c|c|c|}
\hline
\textbf{Variable} & \textbf{Description} & \textbf{Number} \\ \hline
$N_i$ & Initial datapoints & 512 \\ \hline
$N_b$ & Boundary datapoints & 10/segment \\ \hline
$N_r$ & Collocation points & 5000/segment \\ \hline
$N_{iter}$ & Number of ADAM iterations & 10000/segment \\ \hline
\end{tabular}
\caption{Description of training data for Cahn Hilliard equation. 10 time steps/segment have been considered and the amount of collocation points generated remains same and doesn't increase as we progress through time.}
\label{tab:CH_training_data}
\end{table}
Table ~\ref{tab:CH_training_data} describes the values of the  hyper-parameters used in bc-PINN. $N_i$ and $N_b$ refers to the number of points considered to enforce the initial and boundary condition respectively. $N_r$ is the number of collocation points per time segment and $N_{iter}$ is the number of ADAM iterations used to train the neural network per time segment  . \\

\begin{figure}[!htbp]
     \centering
     \subfloat[Solution  (top)  and relative error  (bottom)  via  standard PINN]{\includegraphics[width=1\linewidth,valign=t]{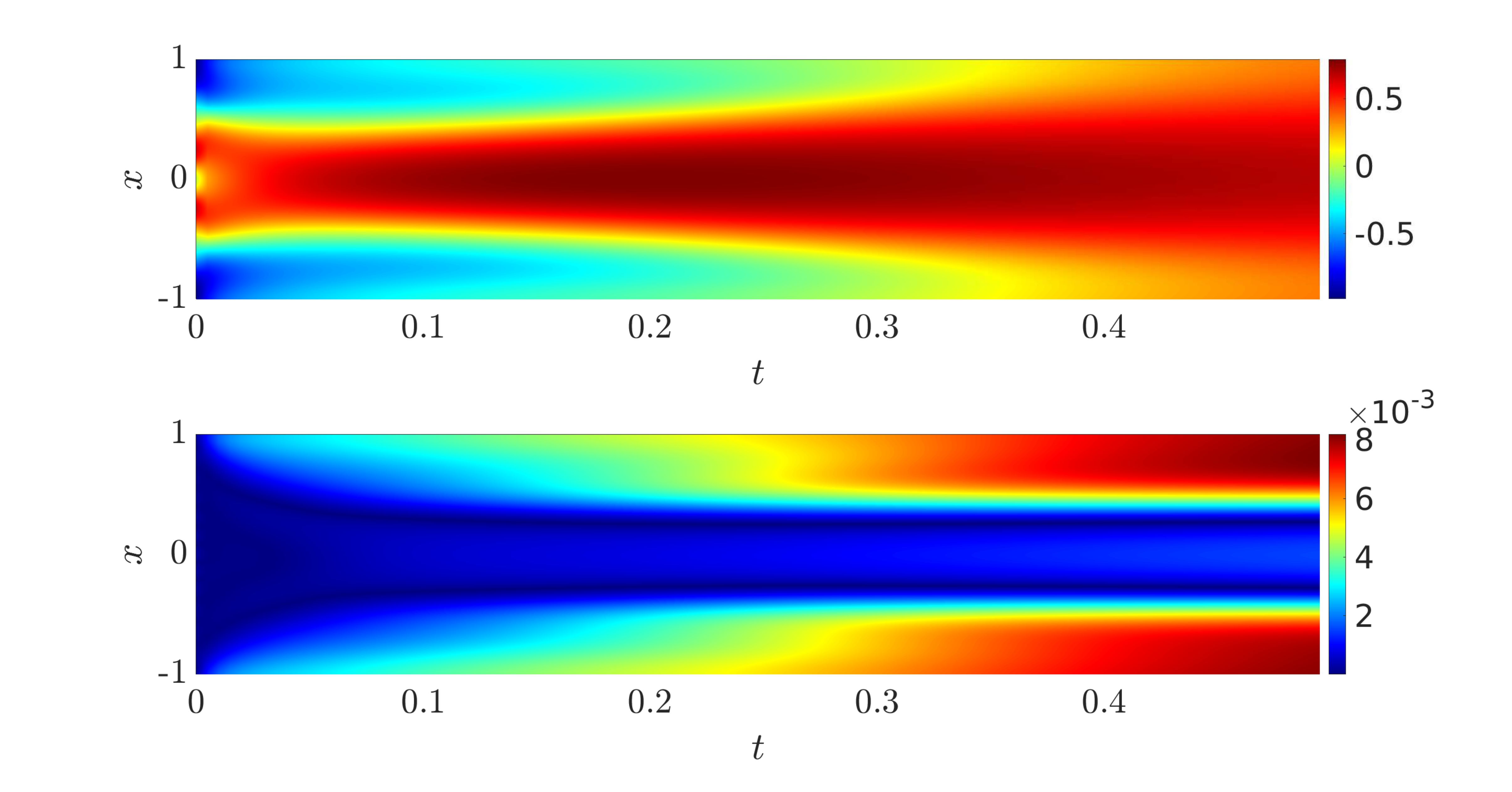}} \quad
     \subfloat[Solution (top) and relative error (bottom) via bc-PINN]{\label{fig:bc-PINN_CH_errplot}\includegraphics[width=1\linewidth,valign=t]{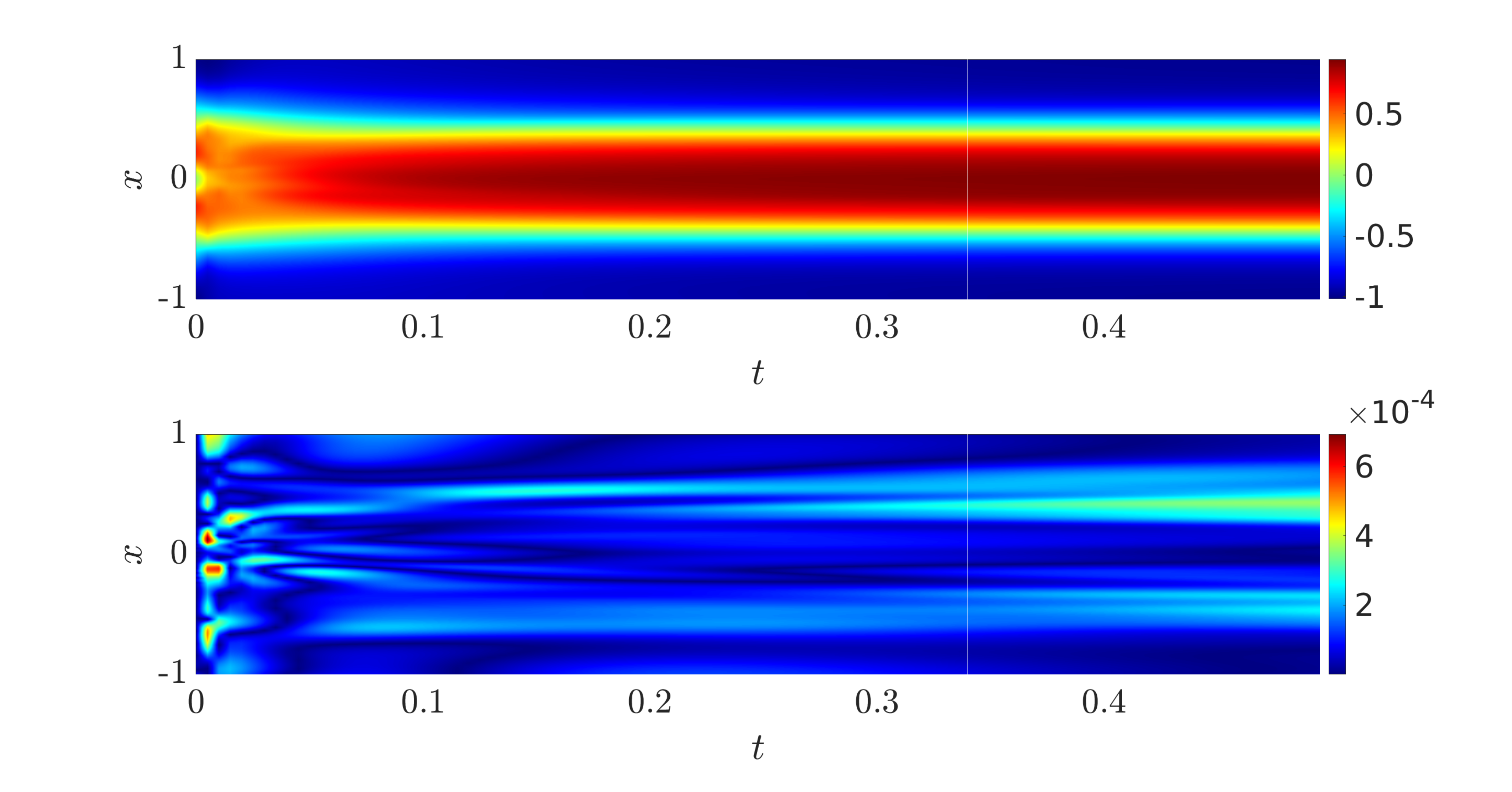}}
     \caption{Solution and error associated with respect to the true solution for Cahn Hilliard equation.} 
     \label{fig:CH_errplot}
\end{figure}

\begin{figure}[!htbp]
     \centering
     \subfloat[Solution (top) and error (bottom) for $\alpha\, \kappa$ = 0.001]{\includegraphics[width=1\linewidth,valign=t]{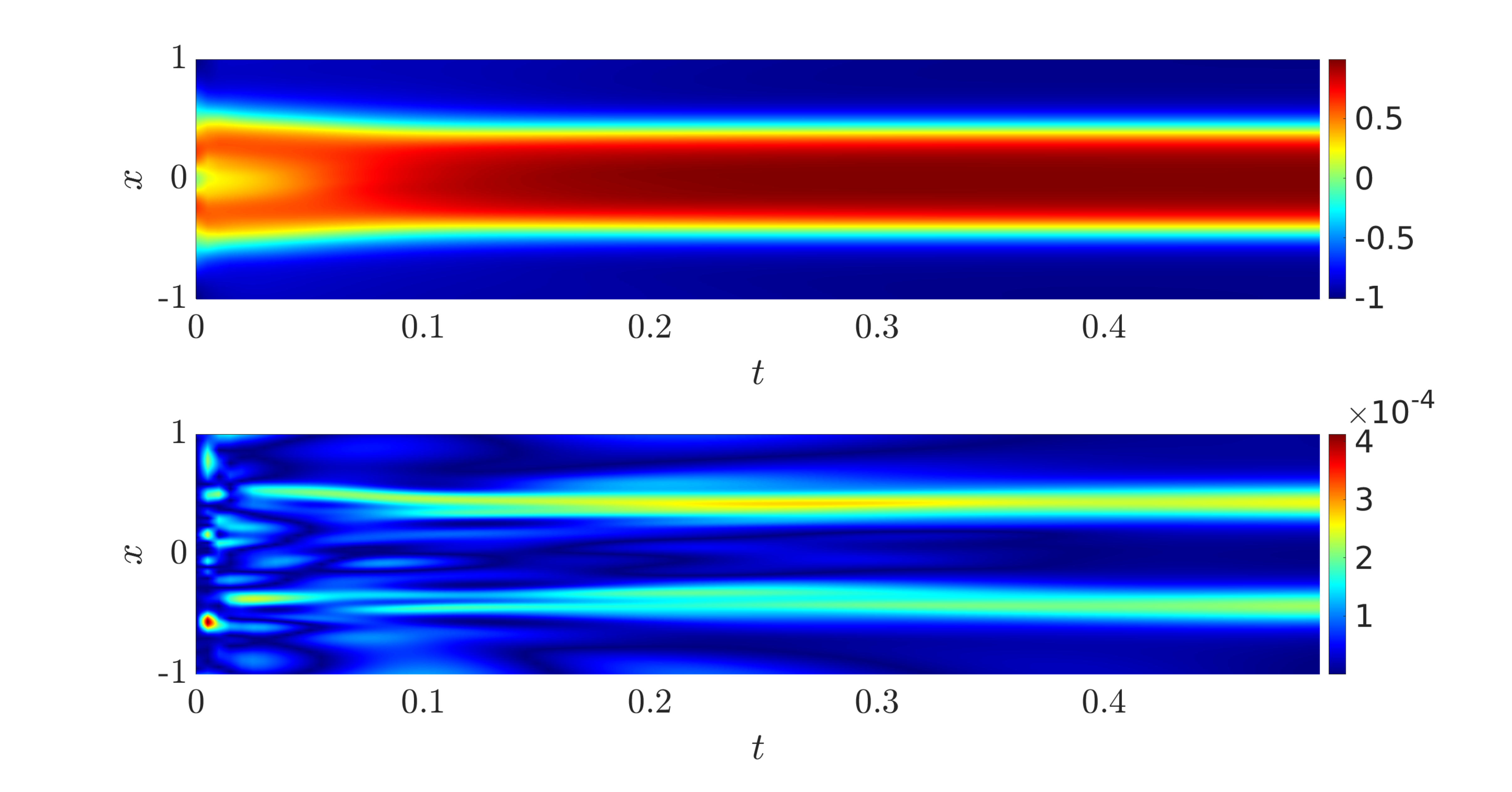}} \quad
     \subfloat[Solution (top) and error (bottom)  $\alpha\, \kappa $ = 0.0005]{\includegraphics[width=1\linewidth,valign=t]{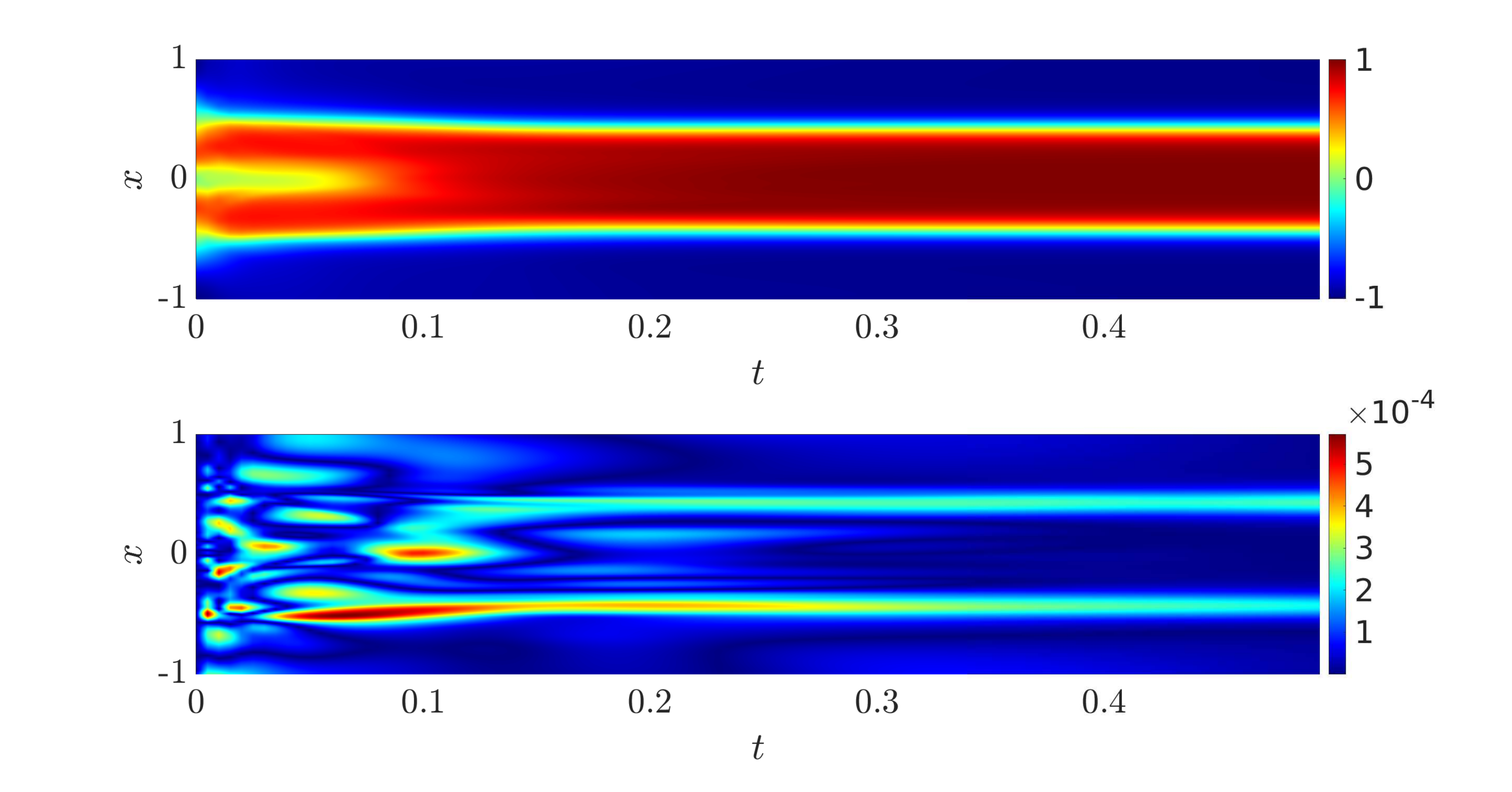}}
    \caption{Solution and relative errors of the Cahn Hilliard equation for different parameters ($\alpha\, \kappa$ of equation~(\ref{eq:21})) obtained by the bc-PINN method.} 
    \label{fig:CH_diff_C1}
\end{figure}

The accuracy of the proposed bc-PINN approach is shown by comparing it against the exact solution obtained by the chebfun method in figure~(\ref{fig:CH_bc-PINN}). This shows that the phase space representation with bc-PINN can closely match the exact solution for the Cahn Hilliard equation. The relative total error ($\varepsilon_{total}$) obtained for the bc-PINN solution is 0.036 whereas for the standard PINN solution the error is 0.8594. It is evident from the error plots given in figure~(\ref{fig:CH_errplot}) that a more accurate solution is obtained by using the bc-PINN compared to standard PINN. The higher accuracy can be accredited to the fact that approximating lower order derivatives using automatic differentiation is much simpler. One key observation to note is that the solution in the $n^{\textrm{th}}$ time segment takes the solution at time $T_{n-1}$  from the $(n-1)^{\textrm{th}}$ time segment as initial condition. Thus,  only the error at the end point in a time segment is propagated to the next time segment. For instance, only the error at the  time $T_{n-1}$ in  $(n-1)^{\textrm{th}}$ time segment is propagated to the next time segment. Errors at all other time steps in $(n-1)^{\textrm{th}}$ time segment does not propagate to the $n^{\textrm{th}}$ time segment. This can be observed in figure~(\ref{fig:bc-PINN_CH_errplot}), even though the error at time $0.01$ is quite high but since this is not the end point of the time segment $[0,0.05]$ it does not propagate with time. The error in the first time segment can be further reduced by using more iterations and the accuracy of the total solution can be improved. To further demonstrate the effectiveness of the current phase space backward compatible training approach, we have taken different values of the parameter ($\alpha \kappa$) and compared the predicted solutions with the exact solutions generated using chebfun which is shown in figure~(\ref{fig:CH_diff_C1}). The proposed phase space representation with bc-PINN approach can be extended to any partial differential equation consisting higher order derivatives and non-linearity. 

 \begin{figure}[!htb]
     \centering
     \subfloat[]{\label{fig:CH_bc-PINN}\includegraphics[width=1\linewidth,valign=t]{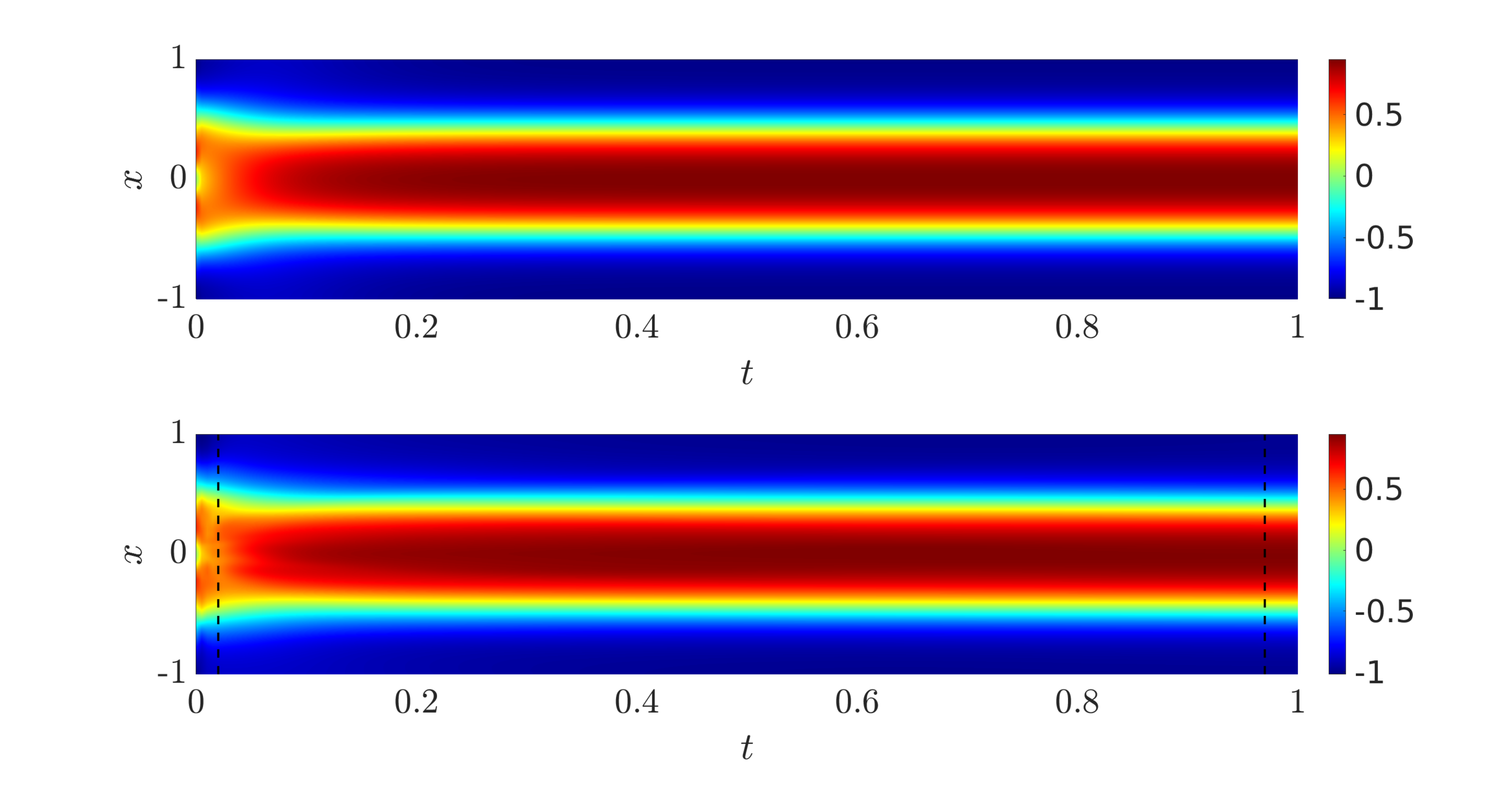}} \quad
     \subfloat[]{\label{fig:CH_5_195}\includegraphics[width=0.8\linewidth,valign=t]{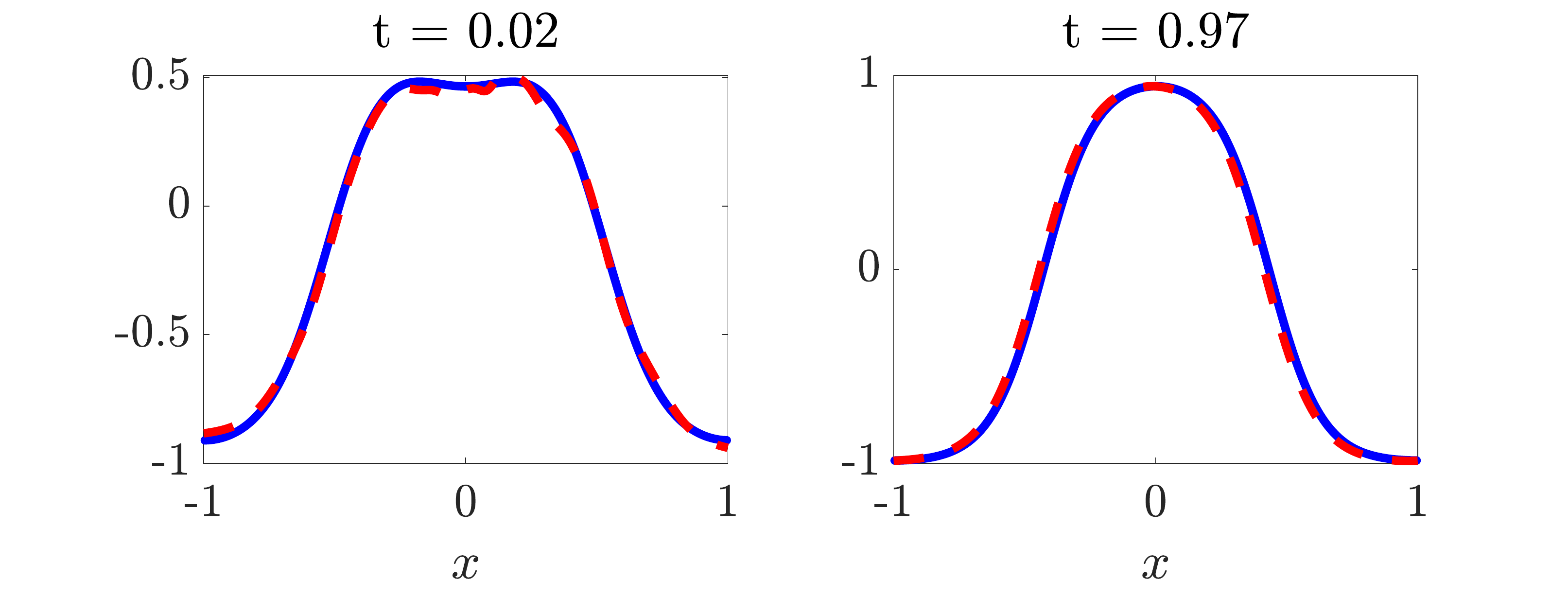}}
     \caption{(a): Exact (Top) and bc-PINN (Bottom) solutions of the Cahn Hilliard equation for the entire spatio--temporal domain. (b): The exact (\blueline) and the  bc-PINN (\redline) solutions at time $t=0.02$ and $t=0.97$.}
     \label{fig:CH_bc-PINN}
\end{figure}


\clearpage
\section{Conclusions}\label{sec:conclusion}
We have proposed a new PINN approach (named as bc-PINN) for solving nonlinear and higher order PDEs. The bc-PINN re-trains the  neural network  over successive time segments while satisfying the solution for all previous time segments. 

Additionally, bc-PINN incorporates different techniques such as  logarithmic residual, phase space representation of the PDE to improve its accuracy. 
The key advantages of bc-PINN are summarized below. \\ 
The proposed bc-PINN method can provide accurate solution for nonlinear and higher--order PDEs such as Cahn Hilliard and Allen Cahn equations, where standard PINN faces difficulties. 

Moreover, the proposed method can achieve high accuracy by using less number of collocation points.
The phase space technique used in the bc-PINN significantly reduces the time required to compute the  derivatives in a higher order PDE. 
Despite the segmentation of the time domain, it requires only one neural network and provides a continuous solution for the entire spatio--temporal domain. The proposed  backward compatibility scheme may  enhance many other machine learning approaches applied to complex systems represented by time dependent PDEs. \\ \\ 

\textbf{Acknowledgments:} 
SG acknowledges the financial support by NSF (CMMI MoMS) under grant number 1937983. We acknowledge Superior, a high-performance computing facility at MTU and Google Colab, a cloud service hosted by Google. This work used the Extreme Science and Engineering Discovery Environment (XSEDE) (allocation number MSS200004), which is supported by the NSF grant number ACI-1548562.

\appendix

\section{Hyper-parameter selection for bc-PINN}\label{sec:HyperParam}
\begin{figure}[!htb]
    \centering
    \includegraphics[scale=0.65]{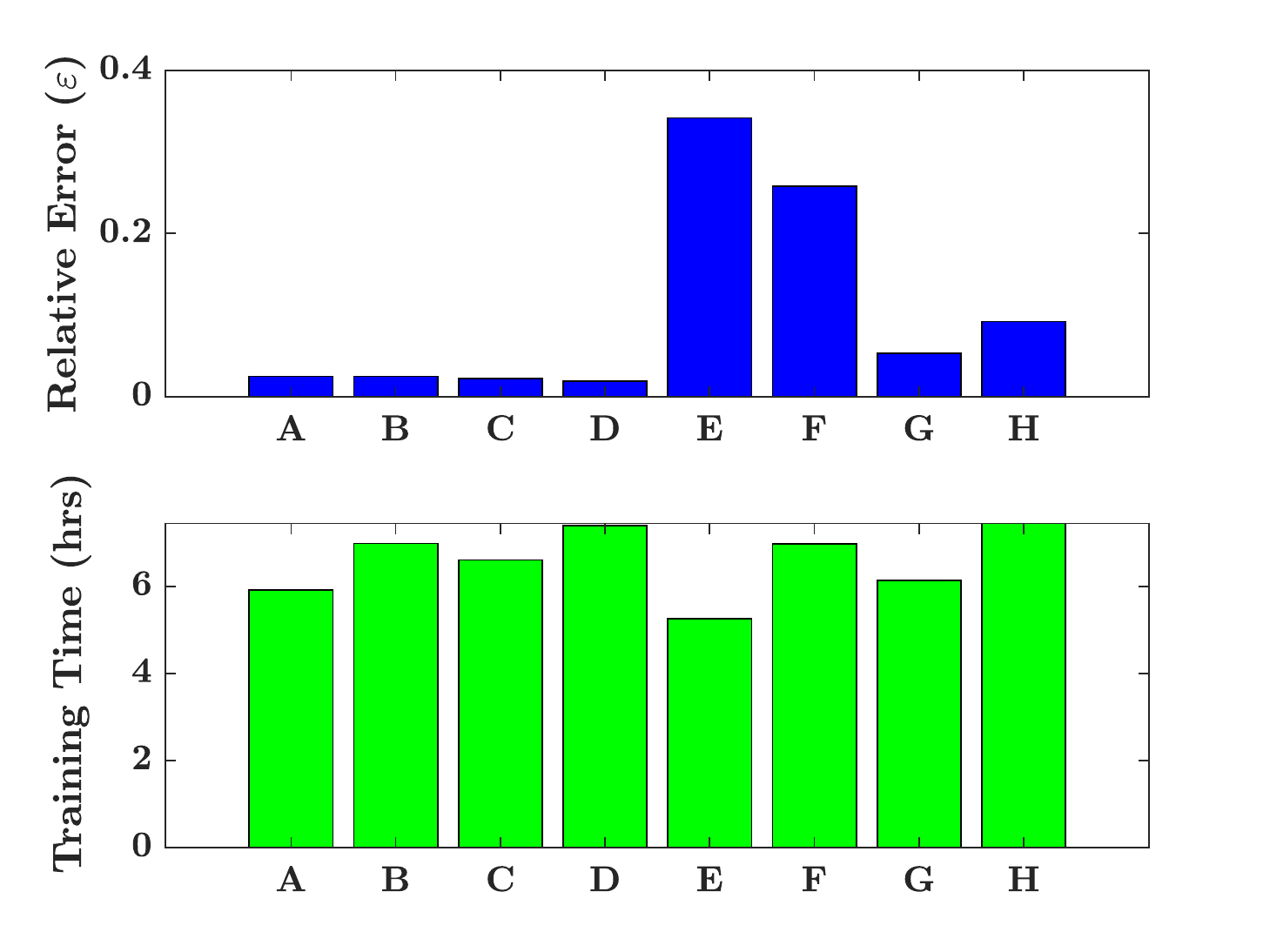}
    \vspace{-15pt}
    \caption{Relative error ($\varepsilon$) and time taken for various models given in table~(\ref{tab:CH_opt_parameters}).} 
    \label{fig:CH_opt_parameters}
\end{figure}

As discussed in section~\ref{sec:AC} \& \ref{sec:CH}, the proposed method there have a number of hyper-parameters like number of ADAM iterations ($N_{iter}$ per segment), time steps per segment, number of collocation points ($N_{r}$) etc. In the current section we choose the Cahn Hilliard equation as the canonical example for all the analysis performed. The accuracy of the bc-PINN's solution depends on proper choice of these hyper-parameters. To optimize each of the hyper-parameters, we considered various cases and metrics like computational time and accuracy. Table~(\ref{tab:CH_opt_parameters}) describes the optimum parameters required to train 100 steps. The optimum parameters are chosen to achieve an accurate solution while balancing the computational cost as shown in figure~(\ref{fig:CH_opt_parameters}). It can be also seen that as the number of collocation points and number of iterations are increased the accuracy increases. So, we choose an segment size of 10 steps with 5000 collocation points and 10000 iterations per time segment as the computational time required in this case is comparably less than the other cases.



 \begin{table}[h!]
 \centering
 \begin{tabular}{|c|c|c|c|}

 \hline
 \textbf{Model} & \textbf{Time steps/segment} & \textbf{$N_r$} & \textbf{$N_{iter}$} \\ \hline
 A & 10 & 5000  & 10000  \\ \hline
 B & 10 & 5000  & 20000  \\ \hline
 C & 10 & 10000 & 10000  \\ \hline
 D & 10 & 10000 & 20000  \\ \hline
 E & 25 & 5000  & 10000  \\ \hline
 F & 25 & 5000  & 20000  \\ \hline
 G & 25 & 10000 & 10000  \\ \hline
 H & 25 & 10000 & 20000  \\ \hline
 \end{tabular}
 \caption{Parameter combinations for choosing the optimum segment size, collocation points and number of ADAM iterations to apply the bc-PINN technique for Cahn Hilliard equation.}
\label{tab:CH_opt_parameters}
\end{table}

\section{bc-PINN with a logarithmic residual for Allen Cahn Equation}\label{sec:bc-PINN_logres}
In this section we show how the bc-PINN with logarithmic residual compares against the standard PINN and bc-PINN without the logarithmic residual. The loss function for the bc-PINN with a logarithmic residual is same as the bc-PINN except the equation(\ref{eq:16}) is replaced by the following loss term for the residual of the PDE:

\begin{align}
    \begin{gathered}
        R := \hat{h}_t - c_{1}^{2}\;\nabla^{2}\hat{h} + f(\hat{h})\\
        \textrm{MSE}_R^{(\ln)} = \,\frac{1}{N_r}\sum_{k=1}^{N_r}\,\ln\left(1 +     (R(x_{k}^{r},t_{k}^{r}))^{2}\right)\,, 
        \qquad (x_{k}^{r},t_{k}^{r}) \in \Omega \times (T_{n-1},T_{n}]
    \label{eq:19}
    \end{gathered}
\end{align}

\begin{figure}[!htb]
     \centering
     \subfloat[bc-PINN]{\label{fig:AC_bcPINN_t_50}\includegraphics[width=0.3\linewidth]{AC_bc-PINN_t_50.pdf}}
     \subfloat[bc-PINN with log residual]{\label{fig:AC_bcPINN_t_50_wlog}\includegraphics[width=0.3\linewidth]{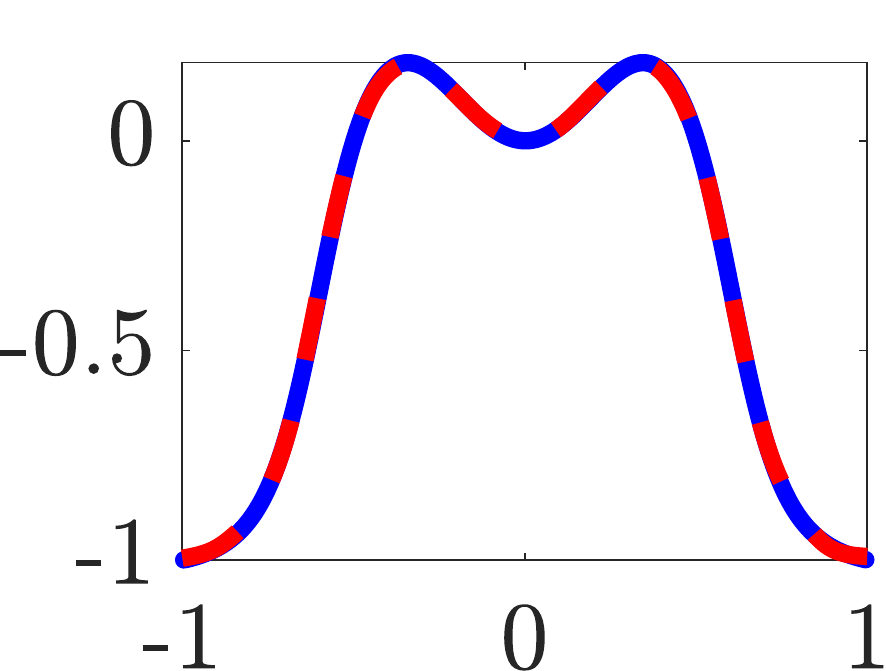}}
     \vspace{10pt}
     \caption{Exact (\blueline) and Predicted (\redline) solution at time t = 0.25}
     \label{fig:AC_bc-PINN-solutions_logres}
\end{figure}

\begin{table}[h!]
\centering
\begin{tabular}{|c|c|} 
\hline
\textbf{Method} & \textbf{Error($\mathbf{\varepsilon}$)} \\ \hline
Standard PINN & 0.9919  \\ \hline
bc-PINN & 0.0701 \\ \hline
bc-PINN with logresidual & 0.03  \\ \hline
\end{tabular}
\caption{\small{Relative errors (equation~(\ref{eq:7})) over the entire domain with respect to Chebfun solution for different methods.}}
\label{tab:AC_error_wlogres}
\end{table}

It turns out that the bc-PINN with a logarithmic residual is more accurate than the bc-PINN.  A possible explanation is that the logarithmic function reduces the relative weight on the $\mathrm{MSE}_R$, which would have larger inaccuracy due to its derivative and nonlinear terms. Thus the initial and boundary terms are satisfied more accurately, which in turn yields a more accurate solution. This explanation is substantiated by the fact that when the logarithmic function is used on all of the four terms in the loss function then the accuracy decreases.

\section{Minimization of the bc-PINN loss function }\label{sec:LossPlots}
\begin{figure}[h!]
    \centering
    \includegraphics[scale=0.8]{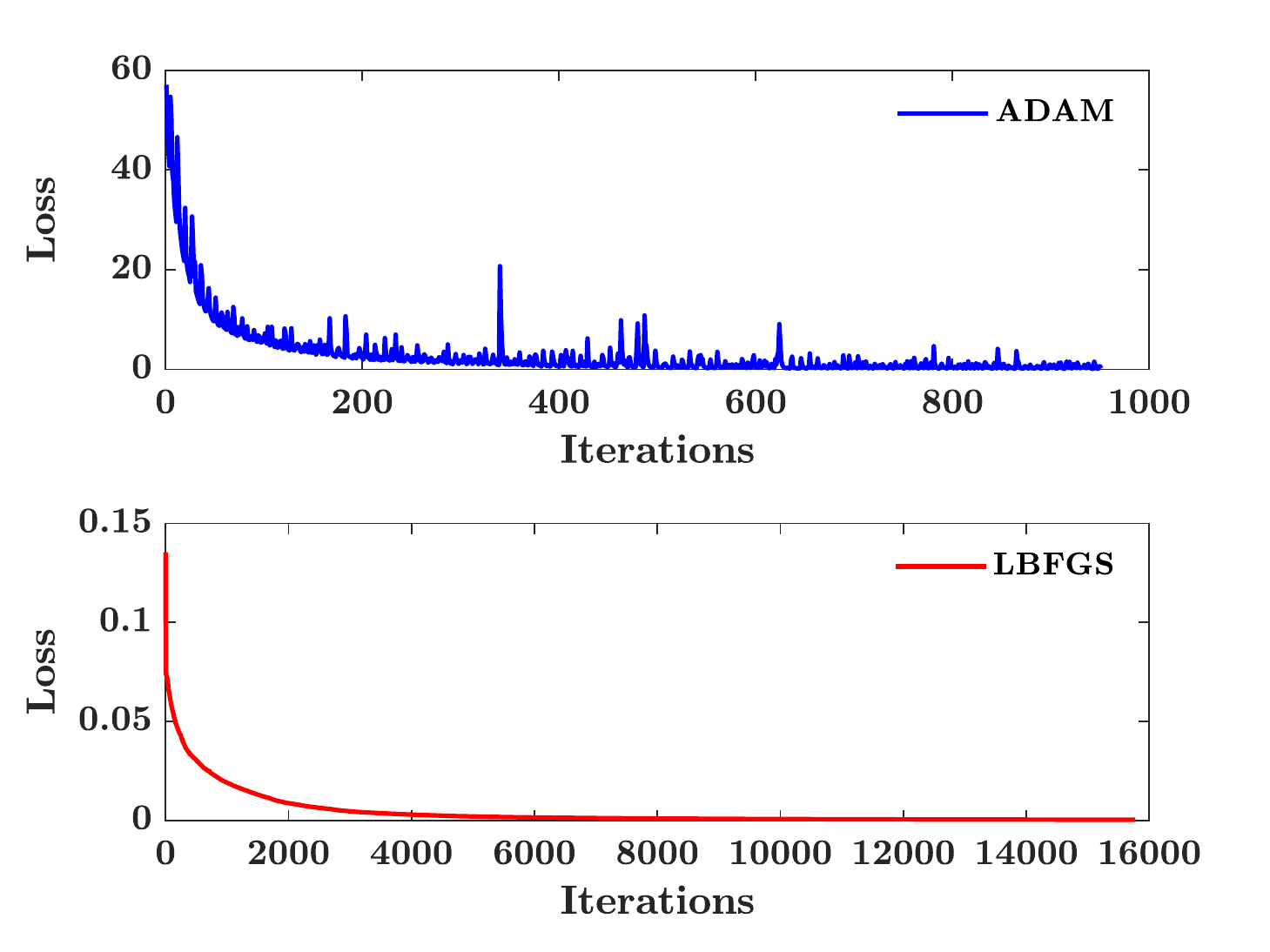}
    \caption{Loss vs Iterations (Top): using ADAM optimizer (Bottom): using the LBFGS optimizer for training the time segment $[0.45,0.5]$ of the Cahn Hilliard equation.}
    \label{fig:CH_loss_90_100}
\end{figure}

In section~\ref{sec:NN} we have described we have mentioned about the learning rates and stopping criteria for the ADAM and LBFGS optimizer are utilized to train the bc-PINN. Here, the minimization of loss function (equation~(\ref{eq:26})) for training the bc-PINN in time segment $[0.45,0.5]$ is given in figure~(\ref{fig:CH_loss_90_100}). 

\newpage
\bibliography{mybibfile}
\end{document}